\documentclass[12pt]{article}

\usepackage{a4wide}
\usepackage{amssymb}
\usepackage{amsfonts}
\usepackage{amsmath}
\input xy
\xyoption{arrow} \xyoption{matrix}

\date{}

\newtheorem{proposition}{Proposition}[section]
\newtheorem{theorem}[proposition]{Theorem}
\newtheorem{lemma}[proposition]{Lemma}

\newtheorem{corollary}[proposition]{Corollary}

\def\der{\partial }

\def\nFM0{{\nu }_{F,M_0}}
\def\nFN0{{\nu }_{F,N_0}}
\def\nGN0{{\nu }_{G,N_0}}

\def\N0{ {\bf N}_0 }

\def\t{\otimes}
\def\g{\gamma}

\def\ra{\rightarrow}

\def\Xpm{X^{\pm }}

\def\s{\sigma}

\def\l1{{\lambda}_1}

\def\a{\alpha}
\def\a0{ {\alpha }_0}
\def\a1{ {\alpha }_1}

\def\l{\lambda}

\def\nFGM0{{\nu }_{F,G,M_0}}

\def\nFN0{{\nu}_{F,N_0}}

\def\sm{{\sigma}^m}

\def\sm1{{\sigma}^{-1}}

\def\smtp1{{\sigma}^{-t+1}}

\def\S1{S^{-1}}

\def\Xpm1{X^{\pm 1}_1}

\def\sPM1{{\sigma }^{\pm 1}}
\def\sMP1{{\sigma }^{\mp 1 }}

\def\d{\delta}

\def\di{{\rm d.ind}}

\def\CA{{\cal A}}

\def\CD{{\cal D}}

\def\Ytm1{Y^{t-1}}
\def\Yim1{Y^{i-1}}

\def\CN{{\cal N}}

\def\k{{\bf k}}

\def\supp{{\rm supp }}

\def\Aut{{\rm Aut}}

\def\Der{{\rm Der }}
\def\ad{{\rm ad }}
\def\dim{{\rm dim }}

\def\ker{ {\rm ker } }

\def\D{ \Delta }
\def\SL2Z{ {\rm SL}_2({\bf Z}) }

\def\Gp1{ G^{1 , 1 } }
\def\P11{ P^{-1 , 1 } }
\def\Pp1{ P^{1 , 1 } }

\def\nCLsr{{}^\nu\kern-2pt {\cal L}^{\sigma , \rho  }}
\def\nP{{}^\nu \kern-2pt P}
\def\nL{{}^\nu\kern-2pt L}
\def\nLL{{}^\nu\kern-2pt \Lambda}
\def\nPsr{{}^\nu\kern-2pt P^{\sigma , \rho  }}
\def\nLsr{{}^\nu\kern-2pt L^{\sigma , \rho  }}
\def\nuCL{{}^\nu\kern-2pt  {\cal L}}
\def\nCLsr{{}^\nu\kern-2pt {\cal L}^{\sigma , \rho  }}
\def\nCL1m{{}^\nu\kern-2pt {\cal L}^{-1 , 1  }}
\def\x1nu{x^\frac{1}{\nu}}
\def\xm1nu{x^{-\frac{1}{\nu}}}

\def\CN{{\cal N}}
\def\ra{\rightarrow }

\def\nAM0{{\nu }_{{\cal A},M_0}}
\def\nAN0{{\nu }_{{\cal A},N_0}}

\def\Der{ {\rm Der }}

\def\det{ {\rm det }}
\def\ad{ {\rm ad }}

\def\ga{\mathfrak{a}}

\def\gn{\mathfrak{n}}
\def\gm{\mathfrak{m}}

\def\derij{\partial_{{\bf i}, {\bf j}}}

\def\di!{\frac{\der^i}{i!}}
\def\dik!{\frac{\der^k_i}{k!}}

\def\dera{\partial^{[\alpha]}}
\def\derb{\partial^{[\beta]}}

\def\xba{x^{[\alpha ]}}

\def\xbi{x^{[i]}}
\def\xbj{x^{[j]}}
\def\ybi{y^{[i]}}
\def\ybj{y^{[j]}}
\def\xbij{x_i^{[j]}}
\def\xbik{x_i^{[k]}}

\def\derik{\der_i^{[k]}}
\def\derij{\der_i^{[j]}}

\def\deril{\der_i^{[l]}}
\def\derjl{\der_j^{[l]}}
\def\derba{\der^{[\alpha ]}}
\def\derbb{\der^{[\beta ]}}
\def\dersba{{\der'}^{[\alpha ]}}

\def\Nn{\mathbb{N}^n}
\def\Nm{\mathbb{N}^m}
\def\Ns{\mathbb{N}^s}
\def\Pnhi{P_{n,\widehat{i}}}
\def\Ahr{A_{\widehat{r}}}
\def\hP{\widehat{P}}
\def\hgm{\widehat{\gm}}
\def\gl{\mathfrak{l}}
\def\NI{\mathbb{N}^{(I)}}
\def\id{{\rm id}}
\begin{document}

\author{V. V. \  Bavula %(chpinv.tex)
}

\title{The inversion formulae for automorphisms of  polynomial algebras and differential operators in prime characteristic}

\maketitle
\begin{abstract}
Let $K$ be an {\em arbitrary} field of characteristic $p>0$, let
$A$ be one of the following algebras: $P_n:= K[x_1, \ldots , x_n]$
is a polynomial algebra, $\CD (P_n)$ is the ring of differential
operators on $P_n$, $\CD (P_n)\t P_m$, the $n$'th {\em Weyl}
algebra $A_n$, the $n$'th {\em Weyl} algebra $A_n\t P_m$ with
polynomial coefficients $P_m$, the power series algebra $K[[x_1,
\ldots , x_n]]$, $T_{k_1, \ldots , k_n}$ is the   subalgebra of
$\CD (P_n)$ generated by $P_n$ and the higher derivations
$\der_i^{[j]}$, $0\leq j <p^{k_i}$, $i=1, \ldots , n$ (where $k_1,
\ldots , k_n\in \mathbb{N}$), $T_{k_1, \ldots , k_n}\t P_m$, an
 {\em arbitrary  central simple} (countably generated) algebra over an
{\em arbitrary} field. {\em The inversion formula} for
automorphisms of the algebra $A$ is found {\em explicitly}.


$${\bf Contents}$$
\begin{enumerate}
\item Introduction. \item Locally nilpotent derivations and their
nil algebras. \item Commuting locally nilpotent derivations with
iterative descents of maximal length. \item A formula for the
inverse of an automorphism that preserves the ring of invariants.
\item  An inversion formula for an automorphism of a central
simple algebra. \item The inversion formula for an automorphism of
a polynomial algebra.\item The inversion formula for an
automorphism $\s \in \Aut_K(\CD (K[x_1, \ldots , x_n]))$. \item
The inversion formula for  $\s \in \Aut_K(\CD (K[x_1, \ldots ,
x_n])\t K[y_1, \ldots , y_m])$. \item The inversion formula for an
automorphism of the $n$'th Weyl algebra $A_n$ (and of $A_n\t
K[y_1, \ldots , y_m]$). \item The inversion formula for $\s \in
\Aut_{K,c}(K[[x_1, \ldots , x_n]])$. \item The inversion formula
for $\s \in \Aut_K(T_{k_1, \ldots , k_n}\t P_m)$.
\end{enumerate}

\end{abstract}

%%%%%%%%%%%%%%%%%% SECTION 1 %%%%%%%%%%%%%%%%%%%%%%%%

\section{Introduction}

Let $K$ be an arbitrary field, $P_n:= K[x_1, \ldots , x_n]$ be a
polynomial algebra over $K$, and $\hP_n:= K[[x_1, \ldots , x_n]]$
be an algebra of formal power series over $K$.

In characteristic {\em zero}, there are several {\em different}
inversion formulae for $\s \in \Aut_K(P_n)$ (\cite{BCW},
\cite{McKayWang88Inv}, \cite{Adj-vdEssen90},\cite{vdEssenPA},
\cite{invfor05}), they are based on different ideas. Besides
applications the interest in inversion formulae stems mainly from
three different sources: {\em in Algebra} - to solve the {\em
Jacobian Conjecture}; {\em in Differential Equations} - to study
solutions of differential equations (and their dependence on
parameters) where on various stages of finding solutions new
coordinates (substitutions) are used; {\em in Analysis} - any set
of `good' functions $ x_1' ,\ldots , x_n'$ in variables $x_1,
\ldots , x_n$ with nonzero {\em Jacobian}, $ {\rm det} (\frac{\der
x_i'}{\der x_j})\neq 0$, defines (locally) coordinates. If, in
addition, the functions $ x_1' ,\ldots , x_n'$ depend on a set of
parameters $\l$ then properties and behaviour of the functions
$$ x_1= x_1(x_1', \ldots , x_n'; \l ),\ldots , x_n= x_n(x_1', \ldots , x_n'; \l
)$$ on $\l $ are of great importance (and difficult to study). A
`nice' inversion formula can help greatly with that sort of
questions. One should also mention the {\em Number Theory}
especially its connections with {\em differential operators} where
the coefficients of differential operators carry valuable
information about the number theoretic question they have
connection with - the inversion formulae in \cite{invfor05} and in
the present paper are given predominantly via differential
operators (or/and powers of derivations in the noncommutative
setting).

In characteristic {\em zero}, the ring of differential operators
$\CD (P_n)$ on the polynomial algebra $P_n$ is  canonically
isomorphic to the Weyl algebra $A_n$ (this is {\em not} the case
in {\em prime} characteristic).  The counter part to the {\em
Jacobian Conjecture} for polynomial algebras    is the {\em
Problem of Dixmier} for the Weyl algebras $A_n$ which asks {\em
whether an algebra endomorphism of $A_n$ is an automorphism}
(\cite{Dix}, Problem 1). Recent results show that these two
problems are essentially the same ($DP_n$ $ \Rightarrow$ $JC_n $,
\cite{BCW}; $JC_{2n}$ $ \Rightarrow$ $DP_n$, \cite{Tsuchi05} and
\cite{Bel-Kon05JCDP}, see also \cite{JC-DP}).  One can even
amalgamate $DP_n$ and $JC_m$ into a single question about algebra
endomorphisms of the algebra $A_n\t P_m$, \cite{invfor05}, though
this question is equivalent to $DP_n+JC_m$, \cite{invfor05}.  The
inversion formula for $\s \in \Aut_K(A_n\t P_m)$ was found in
\cite{invfor05}, it was used to prove the equivalence just
mentioned. In \cite{invfor05}, it was shown that the algebras $\{
A_n\t P_m\}$ are the {\em only} algebras for which the type of
questions like $JC$ or $DP$ makes  sense (i.e. both $JC$ or $DP$
can be reformulated as questions about {\em certain} commuting
locally nilpotent derivations, the algebras $\{ A_n\t P_m\}$ are
the {\em only} algebras that have them).

In characteristic $p>0$, before the present paper I have known no
inversion formula neither  for polynomials (or series) or the Weyl
algebras or the ring of differential operators $\CD (P_n)$. An
 attempt was made by Tsuchimoto \cite{Tsuchi03},
Proposition 1,  to find such a formula for the Weyl algebra $A_n$
using the reduced trace of its division algebra ${\rm Frac} (A_n)$
(in more detail, the `inversion formula' is found for $\s \in
\Aut_K(A_n)$ up to $\s^{-1}|_{Z(A_n)}$ where $Z(A_n)$ is the
centre of the Weyl algebra, it is a polynomial algebra in $2n$
variables).

In the present paper, the inversion formula is found (in which
only algebraic operations are present like  addition or
multiplication) for the algebras mentioned in the Abstract. We
generalized the approach of the paper \cite{invfor05}. The prime
characteristic case is more difficult than the characteristic zero
case and many more situations occur.

The paper is organized as follows: in the first four sections a
necessary machinery is  developed which later is applied in
finding the inversion formulae.

{\bf An idea of finding inversion formula}.  Briefly, {\em it is
to express the identity map via `algebraic' operations} (i.e.
operations which are {\em well-behaved}  when applying an algebra
automorphism): let $A$ be an algebra over a field $K$ and $\s$ be
an algebra automorphism of $A$, let $\{ x^\alpha\}$ be a $K$-basis
of $A$ and suppose that the identity map $\id_A: A\ra A$  can be
written as
$$ \id_A (\cdot) =  \sum \l_{\alpha , y_\alpha }(\cdot )\,
x^\alpha$$ where $\l_{\alpha , y_\alpha }: A\ra K$ are `algebraic'
maps depending `algebraically' on a certain set of elements
$y_\alpha $. Applying $\s$  one has  also the  presentation
$$ \id_A (\cdot) =  \sum \l_{\alpha , \s (y_\alpha ) }(\cdot )\,
\s (x^\alpha ) $$ where $\s (\l_{\alpha , y_\alpha })=\l_{\alpha ,
\s (y_\alpha ) }$ since $\l_{\alpha , \s (y_\alpha ) }$ is an
`algebraic' map. Then applying $\s^{-1}$ and using the fact that
$\l_{\alpha , \s (y_\alpha ) }(A)\subseteq K$  we have the
`formula' for the  inverse map
$$ \s^{-1}(\cdot ) = \sum \l_{\alpha , \s (y_\alpha ) }(\cdot )\,
x^\alpha .$$ To realize this simple idea different algebras
require different means. For {\em central simple} algebras it is
`easy to do it' as the {\em Density Theorem} provides such a
presentation as a limit in the {\em finite topology} of certain
maps given explicitly  (on the other extreme, for {\em
commutative} algebras one has to use {\em differential
operators}).

{\it Example}. Let $M_n(K)$, $n\geq 2$, be the matrix algebra over
 a field $K$, it is a central simple finite dimensional algebra.
 Let $\s\in \Aut_K(M_n(K))$. It it well-known that $\s (x)=
 sxs^{-1}$  for some {\em non-singular}  matrix $s$. Therefore, $\s^{-1} (x)= s^{-1} xs$.  For any
 $x=(x_{ij})\in M_n(K)$,
 $$ x= \sum_{i,j=1}^n x_{ij}e_{ij}= \sum_{i,j=1}^n (x_{ij}E)
 e_{ij}= \sum_{i,j=1}^n (\sum_{k=1}^n e_{ki} x e_{jk})
 e_{ij}$$ where $e_{ij}$ are the matrix units and $E$ is the identity matrix. Equivalently,
 $\id_{M_n(K)} (\cdot )=\sum_{i,j=1}^n (\sum_{k=1}^n e_{ki} (\cdot ) e_{jk})
 e_{ij}.$
 Hence, we have the inversion formula given explicitly
 $$\s^{-1}(\cdot ) =\sum_{i,j=1}^n (\sum_{k=1}^n \s (e_{ki})  (\cdot ) \s (e_{jk}))
 e_{ij}. $$
One can easily verify that this is the inversion formula: for  $\s
(x) = sxs^{-1}$,
\begin{eqnarray*}
&& \sum_{i,j=1}^n (\sum_{k=1}^n \s (e_{ki})  x \s (e_{jk}))
 e_{ij}= \sum_{i,j=1}^n (\sum_{k=1}^n s e_{ki} s^{-1}  x s
e_{jk} s^{-1} )
 e_{ij}=\sum_{i,j=1}^n s(\sum_{k=1}^n  e_{ki} \s^{-1}  (x)
e_{jk}) s^{-1}
 e_{ij}\\
 && = \sum_{i,j=1}^n s (\s^{-1}(x)_{ij}E)
 s^{-1}e_{ij}=\sum_{i,j=1}^n \s^{-1}(x)_{ij}
 e_{ij}= \s^{-1} (x).
\end{eqnarray*}

%%%%%%%%%%%%%%%%%% SECTION 2 %%%%%%%%%%%%%%%%%%%%%%%%

\section{Locally nilpotent derivations and their
nil algebras}

This section is about  the structure of algebras that admit
certain locally nilpotent derivations. These derivations appear
naturally in almost all the algebras we consider in the  paper,
and results of this section  are the key ones in finding various
inversion formulae.

{\bf Iterative $\d$-descents}.  Let $A$ be an algebra over a field
$K$ and let $\d $ be a
 $K$-derivation of the algebra $A$. For any elements $a,b\in A$
 and a natural number $n$, an easy induction argument gives
 $$ \d^n(ab)=\sum_{i=0}^n\, {n\choose i}\d^i(a)\d^{n-i}(b).$$
 It follows that the kernel $A^\d:=\ker \, \d $ of $\d $ is a
 subalgebra (of {\em constants} for $\d $) of $A$ and the union of the
 vector spaces $N:=N(\d ,A)=\cup_{i\geq 0}\, N_i$, $N_i:= \ker \, \d^{i+1}$, is a positively
 {\em filtered} algebra ($N_iN_j\subseteq N_{i+j}$ for
 all $i,j\geq 0$), so-called, the {\em nil algebra} of $\d $.
 Clearly, $N_0= A^\d$ and $N:=\{ a\in A \, | \ \d^n (a)=0$
  for some natural $n=n(a)\}$.

A $K$-derivation $\d $ of
 the algebra $A$ is a {\em locally nilpotent } derivation if for
 each element $a\in A$, $\d^n (a)=0$ for all $n\gg 1$.  A $K$-derivation $\d $ is locally nilpotent iff
 $A=N(\d , A)$.

{\it Definition}. Let $\d$ be a $K$-derivation of an algebra $A$
over an arbitrary  field $K$. A {\em finite} or {\em infinite}
sequence $y=\{ \ybi , 0\leq i\leq l-1\}$ of elements in $A$ where
$y^{[0]}:=1$ is called an {\bf iterative sequence} of {\em length}
$l$ if %\marginpar{it1i}
\begin{equation}\label{it1i}
\ybi \ybj = {i+j \choose i} y^{[i+j]}, \; 0\leq i,j\leq l-1, \;
i+j\leq l-1.
\end{equation}
A sequence $y=\{ \ybi , 0\leq i\leq l-1\}$ is called a $\d$-{\bf
descent} if $y^{[0]}:=1$ and %\marginpar{it1j}
\begin{equation}\label{it1j}
\d (\ybi ) = y^{[i-1]}, \; 0\leq i\leq l-1, \; y^{[-1]}:=0.
\end{equation}
{\it Definition}, \cite{simdernp}. If both conditions (\ref{it1i})
and (\ref{it1j}) hold then the sequence $y$ is called an {\bf
iterative $\d$-descent} of {\em length} $l$.

\begin{lemma}\label{sp2Dec06}%\marginpar{sp2Dec06}
 Let $A$ be an algebra over an arbitrary  field $K$,
$\d $ be a $K$-derivation of $A$.
\begin{enumerate}
\item  If $\{ \xbi , i\geq 0\}$ is a $\d$-descent then $N(\d
,A)=\oplus_{i\geq 0}A^\d \xbi = \oplus_{i\geq 0}\xbi A^\d$ and
$N_n=\oplus_{i=0}^n\, A^\d \xbi =\oplus_{i=0}^n\, \xbi A^\d$ for
all $n\geq 0$. \item If  $\d^l=0$ for some $2\leq l<\infty $, and
$\{ \xbi , 0\leq i <l\}$ is a $\d$-descent of length $l$. Then
$\d^{l-1} \neq 0$ and $N(\d ,A)=\oplus_{i=0}^{l-1}A^\d \xbi =
\oplus_{i=0}^{l-1}\xbi A^\d$ and $N_j=\oplus_{i=0}^j\, A^\d \xbi
=\oplus_{i=0}^j\, \xbi A^\d$ for all $0\leq j <l$.
\end{enumerate}
\end{lemma}

{\it Proof}. $1$. Clearly, $N':= \sum_{i\geq 0}A^\d \xbi\subseteq
N:= N( \d , A)=\cup_{n\geq 0}N_n$ since all $\xbi \in N$ and
$A^\d\subseteq N$. Let us show that the sum in the definition of
$N'$ is a {\em direct} one. Suppose this is not the case, then
there is a nontrivial relation of degree $n>0$,
$$ c_0+c_1x^{[1]}+\cdots +c_nx^{[n]}=0, \;\; c_i\in A^\d ,\;\; c_n\neq 0.$$
We may assume that the degree $n$ of the relation above is the
{\em least} one. Then applying $\d $ to the relation above we
obtain the relation of smaller degree (a contradiction):
$$ c_1+c_2x^{[1]}+\cdots +c_nx^{[n-1]}=0.$$
So, $N'=\oplus_{i\geq 0 } A^\d \xbi$. It remains to prove  that
$N=N'$. It suffices to show that all subspaces $N_i$ belong to
$N'$. We use induction on $i$. The base of the induction is
trivial since $N_0=A^\d$. Suppose that $i>0$, and
$N_{i-1}\subseteq N'$. Let $u$ be an arbitrary element of $N_i$.
Then $\d (u)\in N_{i-1}\subseteq N'$, hence $\d
(u)=\sum_{j=0}^{i-1} c_i\xbj =\d (\sum_{j=0}^{i-1} c_j x^{[j+1]})$
for some $c_i\in A^\d$. Hence, $u-\sum_{j=0}^{i-1} c_j
x^{[j+1]}\in A^\d \subseteq N'$, and so $u\in N'$, which proves
that $N_i\subseteq N'$, and so $N=N'$. It is obvious that
$N_n=\oplus_{i=0}^n\, A^\d \xbi$, $n\geq 0$.

Repeating the argument above for the space $N'':= \sum_{i\geq 0}
\xbi A^\d$ we conclude that $N=N''=\oplus_{i\geq 0}\, \xbi A^\d$
and $N_n=\oplus_{i=0}^n\, \xbi A^\d$, $n\geq 0$.

$2$. $\d^{l-1} (x^{[l-1]})=1\neq 0$ hence $\d^{l-1}\neq 0$. For
the rest, repeat literally the same arguments as in the proof of
statement $1$.  $\Box $

There are elements $a^{ij}_k, b^{ij}_k\in A^\d$ such that
$$ \xbi \xbj = \sum_{k=0}^{i+j}a^{ij}_kx^{[k]}=
\sum_{k=0}^{i+j}x^{[k]}b^{ij}_k, \;\; i,j\geq 0.$$ The elements
 $a^{ij}_k, b^{ij}_k$ are {\em unique} (Lemma \ref{sp2Dec06}). If the field $K$
 has characteristic zero and the element $x:= x^{[1]}$ is fixed
 then the elements $\xbi$, $i\geq 2$ can be chosen as
 $\xbi=\frac{x^i}{i!}$ (in general, the elements $\xbi$ are highly
 non-unique but it is not difficult to describe all
 possibilities: if $\{  {x'}^{[i]}\}$ is another choice of the elements $\{  x^{[i]}\}$ then there
 exist infinite number of scalars $\l_1, \l_2, \ldots $ such that
  $  {x'}^{[i]}=x^{[i]}+\sum_{j=1}^i\l_jx^{[i-j]}$, and vice versa).
  In the characteristic zero case, one can say more about the algebra
 $N(\d , A)$ (see Lemma \ref{dx=1}).

{\bf The projection homomorphisms $\phi$ and $\psi$}.  Given a
ring $R$ and its derivation $d$. The {\em Ore extension} $R[x;d]$
of $R$ is a ring freely generated over $R$ by $x$ subject to the
defining relations: $xr=rx+d(r)$ for all $r\in R$.
$R[x;d]=\oplus_{i\geq 0}Rx^i=\oplus_{i\geq 0}x^iR$ is a left and
right free $R$-module. Given $r\in R$, a derivation $(\ad \,
r)(s):=[r,s]=rs-sr$ of $R$ is called an {\em inner} derivation of
$R$.

\begin{lemma}\label{dx=1}%\marginpar{dx=1}
\cite{invfor05} Let $A$ be an algebra over a field $K$ of
characteristic zero and $\d $ be a $K$-derivation of $A$ such that
$\d (x)=1$ for some $x\in A$. Then $N(\d ,A)=A^\d [x; d]$ is the
Ore extension with coefficients from the algebra $A^\d$, and the
derivation $d$ of the algebra $ A^\d$ is the restriction of the
inner derivation $\ad \, x $ of the algebra $A$ to its subalgebra
$A^\d$. For each $n\geq 0$, $N_n=\oplus_{i=0}^n\, A^\d
x^i=\oplus_{i=0}^n\, x^iA^\d$.
\end{lemma}

If the algebra $A$ is {\em commutative} or the element $x$ is {\em
central} the result above is old and well-known by specialists (in
both cases, the algebra $A$ is a polynomial algebra $A^\d [x]$).

The element $x$ from Lemma \ref{dx=1} yields the iterative
$\d$-descent $\{ \xbi := \frac{x^i}{i!}, i\geq 0\}$ of infinite
length.

\begin{theorem}\label{pw2Dec06}%\marginpar{pw2Dec06}
Let $A$ be an algebra over an arbitrary  field $K$, $\d $ be a
locally nilpotent $K$-derivation of the algebra $A$, and $\{ \xbi
, i\geq 0\}$ be an iterative $\d$-descent.  Then the $K$-linear
maps $\phi :=\sum_{i\geq 0} (-1)^i\xbi \d^i , \psi :=\sum_{i\geq
0} (-1)^i\d^i(\cdot ) \xbi :A\ra A$ satisfy the following
properties.
\begin{enumerate}
\item  The maps $\phi$ and $\psi$ are  homomorphisms of right and
left  $A^\d$-modules respectively. \item The maps $\phi$ and
$\psi$ are projections onto the algebra $A^\d$ (see Lemma
\ref{sp2Dec06}):
\begin{eqnarray*}
\phi : A=A^\d \oplus A_+\ra A^\d \oplus A_+, &\;\; a+b\mapsto a,
\;\; {\rm where}\;\; a\in A^\d, \; b\in A_+:=\oplus_{i\geq 1}\xbi
A^\d ,\\
 \psi : A=A^\d \oplus A_+\ra A^\d \oplus A_+, &\;\; a+b\mapsto a,
\;\; {\rm where}\;\; a\in A^\d, \; b\in A_+:=\oplus_{i\geq 1}A^\d
\xbi .
\end{eqnarray*}
In particular,  ${\rm im} (\phi )={\rm im} (\psi )=A^\d$ and $\phi
(y)=\psi (y)=y$ for all $y\in A^\d$. \item $\phi (\xbi )=\psi
(\xbi )=0$ for all $i\geq 1$. \item For each $a\in A$,
$a=\sum_{i\geq 0 } \xbi \phi (\d^i (a))= \sum_{i\geq 0} \psi
(\d^i(a)) \xbi $.  \item If, in addition, the elements $\{ \xbi
\}$ are central then the maps $\phi$ and $\psi$ are $A^\d$-algebra
homomorphisms.
\end{enumerate}
\end{theorem}

{\it Remark}. If ${\rm char } (K)=0$ this is Theorem  2.2,
\cite{invfor05}.

{\it Proof}. Let us prove the theorem, say, for $\phi$,  for the
map $\psi$ arguments are literally the same with obvious  minor
modifications. The map $\phi$ is well-defined since $\d$ is a
locally nilpotent derivation,  it is  a homomorphism of right
$A^\d$-modules (by the very definition of $\phi$), and $\phi
(y)=y$ for all $y\in A^\d$. For each $j\geq 1$, %\marginpar{px1jb}
\begin{equation}\label{px1jb}
\phi (\xbj )= \sum_{i=0}^j (-1)^i \xbi x^{[j-i]} =
(\sum_{i=0}^j(-1)^i {j\choose i})\xbj = (1-1)^j\xbj
 =0.
\end{equation}
So, the map $\phi$ is a projection onto $A^\d$.

 For each $a= \sum_{i\geq 0} \xbi a_i\in A= \oplus_{i\geq 0} \xbi
 A^\d$  where $a_i\in A^\d$, we have $\phi (\d^i (a))= a_i$, hence
 $a=\sum_{i\geq 0 } \xbi \phi (\d^i (a))$.

If, in addition, the elements $\{ \xbi \}$ are central then the
maps $\phi$ and $\psi$ are $A^\d$-algebra homomorphisms since
$A_+$ is a (two-sided) ideal of the algebra $A$.  $\Box$

The derivation $\d$ from Theorem \ref{pw2Dec06} is locally
nilpotent but {\em not} nilpotent. The next corollary is a similar
result but for a {\em nilpotent} derivation $\d$.

\begin{corollary}\label{1pw2Dec06}%\marginpar{1pw2Dec06}
Let $A$ be an algebra over a  field $K$, $\d $ be a nilpotent
$K$-derivation of the algebra $A$ such that $\d^l=0$ for some
$l\geq 2$, and $\{ \xbi , 0\leq i<l\}$ be an iterative
$\d$-descent. Then the $K$-linear maps $\phi :=\sum_{i= 0}^{l-1}
(-1)^i\xbi \d^i , \psi :=\sum_{i=0}^{l-1} (-1)^i\d^i(\cdot ) \xbi
:A\ra A$ satisfy the following properties.
\begin{enumerate}
\item  The maps $\phi$ and $\psi$ are  homomorphisms of right and
left  $A^\d$-modules respectively. \item The maps $\phi$ and
$\psi$ are projections onto the algebra $A^\d$ (see Lemma
\ref{sp2Dec06}):
\begin{eqnarray*}
\phi : A=A^\d \oplus A_+\ra A^\d \oplus A_+, &\;\; a+b\mapsto a,
\;\; {\rm where}\;\; a\in A^\d, \; b\in A_+:=\oplus_{i=1}^{l-1}
\xbi
A^\d ,\\
 \psi : A=A^\d \oplus A_+\ra A^\d \oplus A_+, &\;\; a+b\mapsto a,
\;\; {\rm where}\;\; a\in A^\d, \; b\in A_+:=\oplus_{i=1}^{l-1}
A^\d \xbi .
\end{eqnarray*}
In particular,  ${\rm im} (\phi )={\rm im} (\psi )=A^\d$ and $\phi
(y)=\psi (y)=y$ for all $y\in A^\d$. \item $\phi (\xbi )=\psi
(\xbi )=0$ for all $i\geq 1$. \item For each $a\in A$,
$a=\sum_{i=0}^{l-1}  \xbi \phi (\d^i (a))= \sum_{i=0}^{l-1}
 \psi (\d^i(a)) \xbi $.  \item If, in addition, the
elements $\{ \xbi \}$ are central  and such that $\xbi \xbj \in
A_+$ for all $i$ and $j$ with $i+j\geq l$ then the maps $\phi$ and
$\psi$ are $A^\d$-algebra homomorphisms.
\end{enumerate}
\end{corollary}

{\it Proof}. Repeat the proof of Theorem \ref{pw2Dec06}. $\Box $

{\bf The ring of differential operators  $\CD (P_n)$ on a
polynomial algebra}. If ${\rm char} (K)=p>0$ then {\em the ring
$\CD (P_n)$ of differential operators on a polynomial algebra}
$P_n:= K[x_1, \ldots , x_n]$ is a $K$-algebra generated by the
elements $x_1, \ldots , x_n$ and {\em commuting higher
derivations} $\derik :=\frac{\der_i^k}{k!}$, $i=1, \ldots , n$ and
$k\geq 1$,  that satisfy the following {\em defining} relations:
%\marginpar{DPndef}
\begin{equation}\label{DPndef}
[x_i,x_j]=[\derik , \derjl ]=0,\;\;\; \derik \deril ={k+l\choose
k}\der_i^{[k+l]}, \;\;\; [\derik , x_j]=\d_{ij}\der_i^{[k-1]},
\end{equation}
 for all
$i,j=1, \ldots , n$ and $k,l\geq 1$ where $\d_{ij}$ is the
Kronecker delta and $\der_i^{[0]}:=1$.
$\der_i^{[1]}=\der_i=\frac{\der}{\der x_i}\in \Der_K(P_n)$,
$i=1,\ldots , n$.  It is convenient to set $\der_i^{[-1]}:=0$. The
algebra $\CD (P_n)$ is a {\em simple} algebra. Note  that the
algebra $\CD (P_n)$ is {\em not} finitely generated and {\em not}
(left or right) Noetherian, it {\em does not} satisfy finitely
many defining relations.
$$\CD (P_n)= \bigoplus_{\alpha , \beta \in \Nn }Kx^\alpha \derbb
=\bigoplus_{\alpha , \beta \in \Nn }K \derbb x^\alpha
  , \;\; {\rm where} \;\; x^\alpha := x_1^{\alpha_1}\cdots
x_n^{\alpha_n}, \; \derbb := \der_1^{[\beta_1]}\cdots
\der_n^{[\beta_n]}.$$ The algebra $\CD (P_n)$ admits the {\em
canonical filtration} $F:=\{ F_i\}_{i\geq 0}$ where
$F_i:=\oplus_{|\alpha |+ |\beta |\leq i }Kx^\alpha \derbb $ $
=\oplus_{ |\alpha |+ |\beta |\leq i}K \derbb  x^\alpha $, $
|\alpha |:= \alpha_1+\cdots + \alpha_n$, $\dim_K(F_i)={i+2n\choose
2n}$ for all $i\geq 0$, $\CD
 (P_n)=\cup_{i\geq 0} F_i$, $F_0=K$ and $F_iF_j\subseteq F_{i+j}$
for all $i,j\geq 0$.

For each $i=1, \ldots , n$, the inner derivations $-\ad \, x_i$ of
the algebra $\CD (P_n)$ and the higher derivations $\{ \derij \}$
satisfy the conditions of Theorem \ref{pw2Dec06} (i.e. the
sequence $\{ \derij , j\geq 0\}$ is an {\em iterative $(-{\rm ad}
\, x_i)$-descent}), and the derivations $\der_i:= \der_i^{[1]}$ of
the algebra $P_n$ satisfy the conditions of Corollary
\ref{1pw2Dec06} ($\der_i^p=0$, $\der_i(x_i)=1$, and $\{
\frac{x_i^j}{j!}, 0\leq j <p\}$ is an iterative $\der_i$-descent
of length $p$).

Let $A$ be an algebra over a field $K$ of characteristic $p>0$.
For each $x\in A$, %\marginpar{adxp}
\begin{equation}\label{adxp}
({\rm ad} \, x)^p= {\rm ad }  (x^p)
\end{equation}
since, for any $ a\in A$, ${\rm ad }  (x^p)(a)= [ x^p,
a]=\sum_{i=1}^p{p\choose i} ({\rm ad}\, x)^i (a)x^{p-i}=({\rm
ad}\, x)^p(a)$.

\begin{theorem}\label{lp2Dec06}%\marginpar{lp2Dec06}
Let $A$ be an algebra over a field $K$of characteristic $p>0$, $\d
$ be a $K$-derivation of the algebra $A$ such that $\d^p=0$ and
$\d (x)=1$ for some element $x\in A$. Then
\begin{enumerate}
\item $A=\oplus_{i=0}^{p-1}A^\d x^i= \oplus_{i=0}^{p-1}x^iA^\d$
and $A=N(\d , A)=\cup_{n=0}^{p-1}N_n$ where $N_n=
\oplus_{i=0}^nA^\d x^i=\oplus_{i=0}^nx^iA^\d$ for $n=0, 1, \ldots
, p-1$.
 \item  $[x, A^\d ]\subseteq  A^\d$.
 \item Let $\CA := A^\d [ t; d]$ be an Ore extension of the
 algebra $A^\d$ where $d:= \ad (x)|_{A^\d}$. Then $t^p-x^p$ is a
 central element of the algebra $\CA$ and the algebra $A$ is
 canonically isomorphic to the factor algebra $\CA / (t^p-x^p)$.
\end{enumerate}
\end{theorem}

{\it Remark}. If, in addition, the algebra $A$ is {\em
commutative}, the first statement of Theorem \ref{lp2Dec06} is
Theorem 27.3, \cite{Ma}.

{\it Proof}. Statement 1 is a particular case of Lemma
\ref{sp2Dec06}.(2) with the iterative $\d$-descent $\{ \xbi :=
\frac{x^i}{i!}, 0\leq i <p\}$ of length $p$.

For each $c\in A^\d $, $\d ([x,c])= [\d (x), c]+[x, \d (c)]=
[1,c]+[x,0]=0$, hence $[x,A^\d ]\subseteq A^\d $ and $d$ is a
derivation of the algebra $A^\d$. This proves statement $2$.

% $\d
%(x^p)=px^{p-1}=0$ implies $x^p\in A^\d$, and so the $K$-subalgebra
%$N'$ of $A$ generated by $A^\d$ and $x$ is equal to
%$\sum_{i=0}^{p-1}A^\d x^i= \sum_{i=0}^{p-1}x^iA^\d$. Let us prove
%that, say, the first sum is direct: given a nontrivial relation of
%degree $1\leq n \leq p-1$, $ c_0+c_1x+\cdots +c_nx^n=0$ for some
%elements $c_i\in A^\d$ such that $c_n\neq 0$. Applying $\d^n$ to
%the relation we have a contradiction: $0=n!c_n\neq 0$. A similar
%argument proves that the second sum is direct.
%Let us prove that $A=N'$. Note that $\d $ is a  nilpotent
%derivation of the algebra $A$ and $\d^p=0$. Therefore,
%$A=\cup_{i=0}^{p-1}N_i$ where $N_i:= \ker (\d^{i+1})$. It is
%sufficient to prove that all $N_i\subseteq N'$. We use induction
%on $i$. The case $i=0$ is obvious, $N_0=A^\d \subseteq N'$.
%Suppose that $0<i<p-1$ and $N_{i-1}\subseteq N'$. Let $u$ be an
%arbitrary element of $N_i$. Then $\d (u)\in N_{i-1}$, and so $\d
%(u)=\sum_{j=0}^{i-1}c_jx^j= \d
%(\sum_{j=0}^{i-1}(j+1)^{-1}c_jx^{j+1})$ for some $c_j\in A^\d$
%which implies that $u-\sum_{j=0}^{i-1}(j+1)^{-1}c_jx^{j+1}\in
%A^\d$, and so $u\in N'$. This proves that $N'=A$.  Clearly, $N_n=
%\oplus_{i=0}^nA^\d x^i=\oplus_{i=0}^nx^iA^\d$ for $n=0, 1, \ldots
%, p-1$. This finishes the proof of statement 1.

It remains to prove  statement $3$. Consider the Ore extension
$\CA := A^\d [t; d]$. There is a natural algebra {\em epimorphism}
$\CA \ra A$, $t\mapsto x$, $c\mapsto c$ for all $c\in A^\d$. The
epimorphism is obviously an $A^\d$-{\em module} epimorphism. Since
both one-sided ideals $\CA (t^p-x^p)$ and $(t^p-x^p)\CA$ of the
algebra $\CA$ belong to the kernel of the epimorphism and,
obviously, there are isomorphisms $\CA /\CA (t^p-x^p)\simeq A$ and
$\CA /(t^p-x^p)\CA \simeq A$ of left and right $A^\d$-modules
respectively (induced by the epimorphism $\CA \ra A$), we must
have the equality $\CA (t^p-x^p)=(t^p-x^p)\CA$. Since, for each
$c\in A^\d$, the degree in $t$ of the commutator $[ t^p-x^p, c]\in
\CA$ is {\em strictly less} than $p$ we must have $[ t^p-x^p,
c]=0$ as the commutator belongs to the kernel of the epimorphism
above and, by statement $1$, $A= \oplus_{i=0}^{p-1} A^\d x^i$. One
can prove this fact also directly using (\ref{adxp}) repeatedly:
$$ [ t^p-x^p, c]= ({\rm ad }\, t)^p (c) - [ x^p, c]= ({\rm ad }\, x)^p (c) - [ x^p, c]= [ x^p, c]- [ x^p,
c]=0.$$
 Now,
$[t^p-x^p, t]=[t,x^p]=d(x^p)=[x,x^p]=0$. This means that the
element $t^p-x^p$ belongs to the centre of the algebra $\CA$
(since it commutes with generators of $\CA$) and $A\simeq \CA
/(t^p-x^p)$. $\Box $

The next corollary describes the algebras from Lemma
\ref{lp2Dec06}.

\begin{corollary}\label{clp2Dec06}%\marginpar{clp2Dec06}
Let $A$ be an algebra over a field $K$ of characteristic $p>0$.
Then the following statements are equivalent.
\begin{enumerate}
\item There exists a $K$-derivation $\d $ of the algebra $A$ such
that $\d^p=0$ and $\d (x)=1$ for some element $x\in A$. \item The
algebra $A$ is isomorphic to the factor algebra $\CA / \CA
(x^p-\alpha )$ of an Ore extension $\CA =B[x ; d]$ at the central
element $x^p-\alpha$ for some $\alpha \in B^d$ where $d\in
\Der_K(B)$ such that $d^p=\ad (\alpha )$ (more precisely, for any
choice of $\alpha \in B^d$ and a derivation $d$ of $B$ such that
$d^p = {\rm ad } (\alpha )$, the element $x^p-\alpha$ is a central
element of $\CA$). In this case, the derivation $\d$ (from
statement 1) may be chosen as follows: $\d (B)=0$ and $\d (x)=1$
(then $\d^p=0$).
\end{enumerate}
\end{corollary}

{\it Proof}. $(1\Rightarrow 2)$ Theorem  \ref{lp2Dec06}.

$(2\Rightarrow 1)$ First, let us prove that the element
$x^p-\alpha $ is central. For any $b\in B$,
$[x^p,b]=d^p(b)=[\alpha , b]$ since $ d^p= {\rm ad} (\alpha )$,
and so $[x^p-\alpha , b]=0$. Finally, $[x^p-\alpha ,x]=d(\alpha
)=0$ since $\alpha\in \ker \, d$. So, the element $x^p-\alpha$
commutes with the generators $B$ and $x$ of the algebra $A$.
Therefore, it is a central element of the algebra $A$.

Consider the  derivation $\d \in \Der_B(A)$ given by the rule $\d
(x)=1$. It is well-defined since, for each $b\in B$, $\d ([x, b])
= [1, b]+[x, \d (b)]=0= \d ( d(b))$ and $\d
(x^p-\alpha)=px^{\alpha -1}-\d (\alpha )=0$. Since $\d^p(x^i)=0$
for all $i=1, \ldots , p-1$, we must have $\d^p=0$. $\Box $

%%%%%%%%%%%%%%%%%% SECTION 3 %%%%%%%%%%%%%%%%%%%%%%%%

\section{Commuting locally nilpotent derivations with
iterative descents of maximal length}

{\it Definition}. Let $\d$ be a nonzero locally nilpotent
$K$-derivation of an algebra $A$ over an arbitrary field $K$ and
$y:= \{ \ybi \}$ be an iterative $\d$-descent. We say that $y$ is
{\em of maximal length} (or {\em has maximal length}) if $y$ has
{\em infinite} length in the case when $\d$ is {\em not} nilpotent
derivation, and $y$ has length $l$ if $\d^l =0$ and $ \d^{l-1}
\neq 0$.

Let $A$ be an algebra over a field $K$, $\d := \{ \d_i, i\in I\}$
be a {\em non-empty} set of {\em commuting, locally nilpotent}
$K$-derivations of the algebra $A$ such that

$(i)$ for each $a\in A$, $\d_i(a)=0$ for {\em almost all} $i\in I$
(i.e. all but finitely many $i$), and

$(ii)$ for each $i\in I$, there is an iterative $\d_i$-descent
$x_i:=\{ \xbij \}$ of {\em maximal length}, say $l_i$,  such that
$ \{ \xbij \}\subseteq \cap_{i\neq k\in I} A^{\d_k}$.

{\it Remark}. If the set $I$ is {\em finite} then the condition
$(i)$ is vacuous.

 The set $\d$ is a disjoint union of two subsets
$\d = \gn \cup \gl$ where the set $\gn $ contains all the {\em
nilpotent} derivations and the set $\gl$ contains all the {\em
non-nilpotent} derivations. It is possible that one of this sets
is an empty set.

Let $\mathbb{N}=\{ 0,1, 2, \ldots \}$ be the monoid of natural
numbers, $\NI$ be the direct sum of $I$'th copies of the monoid
$\mathbb{N}$,
$$E= E(\d ):= \{ \alpha = (\alpha_i)\in \NI \, | \,
0\leq \alpha_i\leq l_i-1\; {\rm  where}\;  l_i\;{\rm  is \; the \;
length\;  of}\;  x_i\},$$
 it is the set of all `possible exponents' for $\xba$ (see below). Without
loss of generality we may assume that the set $I$ is a {\em
well-ordered} set with respect to an ordering $<$. For each
$\alpha = (\alpha_i )\in E$, consider the {\em ordered} product,
$$ \xba := \begin{cases}
1& \text{if $\alpha =0$},\\
x_{i_1}^{[\alpha_{i_1}]}\cdots x_{i_n}^{[\alpha_{i_n}]} & \text{if $\alpha \neq 0$ and}\\
\end{cases}$$
$i_1<\cdots <i_n$,   and $\alpha_{i_1}, \ldots
 ,\alpha_{i_n}$ are the only {\em nonzero} coordinates of $\alpha$,
the set $\{ \alpha_{i_1}, \ldots , \alpha_{i_n}\}$ is called the
{\em support} of $\alpha $.

For each $i\in I$, consider the maps $\phi_i, \psi_i :A\ra A$ from
Theorem \ref{pw2Dec06} and Corollary \ref{1pw2Dec06}:

(a) if $\{ \xbij , j\geq 0\}$ is an {\em infinite} $\d_i$-descent
then
$$ \phi_i := \sum_{j\geq 0} (-1)^j \xbij \d_i^j\;\; {\rm and}\;\; \psi_i := \sum_{j\geq 0} (-1)^j
\d_i^j(\cdot ) \xbij,$$

(b) if $\{ \xbij , 0\leq j\leq l_i-1\}$ is a {\em finite}
$\d_i$-descent then
$$ \phi_i := \sum_{j=0}^{l_i-1} (-1)^j \xbij \d_i^j\;\; {\rm and}\;\; \psi_i := \sum_{j=0}^{l_i-1} (-1)^j
\d_i^j(\cdot ) \xbij.$$

\begin{theorem}\label{9Apr06}%\marginpar{9Apr06}
Let $A$ be an algebra over a  field $K$, $(I, <)$ be a non-empty
well-ordered set,  $\d := \{ \d_i, i\in I\} $ be a set of
commuting locally nilpotent $K$-derivation of the algebra $A$ such
that, for each $a\in A$, $\d_i(a)=0$ for almost all $i\in I$.
Suppose that, for each $i\in I$, there exists an iterative
$\d_i$-descent $\{ \xbij \}$ of maximal length such that $\{ \xbij
\} \subseteq \cap_{i\neq k\in I} A^{\d_k}$. Then
\begin{enumerate}
\item $A= \oplus_{\alpha \in E} \xba A^\d= \oplus_{\alpha \in E}
A^\d \xba$ where $A^\d := \cap_{i\in I} A^{\d_i}$. \item Consider
the  ordered products of maps  %\marginpar{psphb}
\begin{equation}\label{psphb}
\phi := \prod_{i\in I} \phi_i, \;  \psi := \prod_{i\in I} \psi_i :
A\ra A,
\end{equation}
given by the rules:  $\phi (1):=1$, $\psi (1):=1$, and, for each
$0\neq \alpha \in E$, with $\supp (\alpha ) =\{ \alpha_{i_1},
\ldots , \alpha_{i_n}\}$ where $i_1<\cdots <i_n$,
$$ \phi (\xba ) := \phi_{i_n} \cdots \phi_{i_1} (
x_{i_1}^{[\alpha_{i_1}]}\cdots x_{i_n}^{[\alpha_{i_n}]}), \;\;
\psi (\xba ) := \psi_{i_1} \cdots \psi_{i_n} (
x_{i_1}^{[\alpha_{i_1}]}\cdots x_{i_n}^{[\alpha_{i_n}]}).$$ Then
the maps $\phi$ and $\psi$ are  homomorphisms of right and left
$A^\d$-modules respectively. \item The maps $\phi$ and $\psi$ are
projections onto the algebra $A^\d$:
\begin{eqnarray*}
\phi : A=A^\d \oplus A_+\ra A^\d \oplus A_+, &\;\; a+b\mapsto a,
\;\; {\rm where}\;\; a\in A^\d, \; b\in A_+:=\oplus_{0\neq \alpha
\in E}\xba
A^\d ,\\
 \psi : A=A^\d \oplus A_+\ra A^\d \oplus A_+, &\;\; a+b\mapsto a,
\;\; {\rm where}\;\; a\in A^\d, \; b\in A_+:=\oplus_{0\neq \alpha
\in E}A^\d \xba .
\end{eqnarray*}
In particular,  ${\rm im} (\phi )={\rm im} (\psi )=A^\d$ and $\phi
(y)=\psi (y)=y$ for all $y\in A^\d$. \item $\phi (\xba )=\psi
(\xba )=0$ for all $0\neq \alpha \in E$. \item For each $a\in A$,
$a=\sum_{\alpha \in E } \xba \phi (\d^\alpha (a))= \sum_{\alpha
\in E} \psi (\d^\alpha (a)) \xba $ where $\d^\alpha := \prod_{i\in
I}  \d_i^{\alpha_i}$, a finite product. \item If, in addition, all
the elements $\{ \xbij \}$ are central and for each finite length
sequence $\{ \xbij , 0\leq j<l_i\}$: $\xbij \xbik \in A_+$ for all
$j$ and $k$
 with $j+k\geq l_i$, then the maps $\phi$ and $\psi$ are $A^\d$-algebra homomorphisms.
\end{enumerate}
\end{theorem}

{ {\it Proof}. We have obvious symmetry between $\phi $ and
$\psi$, and between statements about left $A^\d$-modules and right
 $A^\d$-modules. Let us, say, prove the statements about the map
$\phi$ and about right $A^\d$-modules.

$1$. The sum $A':= \sum_{\alpha \in E} \xba A^\d$ is a direct sum,
$A':= \oplus_{\alpha \in E} \xba A^\d$, as follows at once from
the fact that $\d^\alpha (\xba )=1$ for all $\alpha \in E$. We
have to prove that $A'=A$. Let $a\in A$. We have to show that
$a\in A'$. If $a\in A^\d\subseteq A'$ then there is nothing to
prove. So, let $a\not \in A^\d$. By the assumption, there are only
finitely many derivations, say $\d_1, \ldots , \d_n$, such that
$\d_i(a)\neq 0$ (to save on notation  we may assume that $\{
1<\cdots <n\} \subseteq  I$). Applying step by step either Theorem
\ref{pw2Dec06} or Corollary \ref{1pw2Dec06}, we have (where
$A^{\d_1, \ldots , \d_k} := \cap_{j=1}^k A^{\d_j}$)
%\marginpar{An1s}
\begin{equation}\label{An1s}
A= \bigoplus_{i_1} x_1^{[i_1]}A^{\d_1}=\bigoplus_{i_1, i_2}
x_1^{[i_1]}x_2^{[i_2]}A^{\d_1, \d_2}=\cdots =
\bigoplus_{i_1,\ldots ,  i_n} x_1^{[i_1]}\cdots
x_n^{[i_n]}A^{\d_1, \ldots , \d_n},
\end{equation}
we have used the facts that the derivations $\d_1, \ldots , \d_n$
{\em commute} and that $\d_i (x_j^{[k]})=0$ if $i\neq j$. The
element $a$ is a unique finite sum $a= \sum x_1^{[i_1]}\cdots
x_n^{[i_n]}\l_{i_1, \ldots , i_n}$ for some $\l_{i_1, \ldots ,
i_n}\in A^{\d_1, \ldots , \d_n}$. For each $k\in I\backslash \{ 1,
\ldots , n\}$,
$$ 0=\d_k (a)=\sum x_1^{[i_1]}\cdots x_n^{[i_n]}\d_k(\l_{i_1, \ldots , i_n}), $$
therefore $\d_k(\l_{i_1, \ldots , i_n})=0$ by (\ref{An1s}), i.e.
all $\l_{i_1, \ldots , i_n}\in A^\d$. This proves the equality $A=
A'$.

$2-4$. By the very definition, the map $\phi$ is a homomorphism of
right $A^\d$-modules, $\phi (1):=1$, and, for each $0\neq \alpha
\in E$ with $\supp (\alpha ) =\{ \alpha_{i_1} , \ldots ,
\alpha_{i_n}\}$ and $i_1<\cdots <i_n$,
$$ \phi (\xba )= \phi_{i_1}(x_{i_1}^{[\alpha_{i_1}]})\cdots
\phi_{i_n}(x_{i_n}^{[\alpha_{i_n}]})=\d_{0, \alpha_{i_1}}\cdots
\d_{0, \alpha_{i_n}}= \d_{0, \alpha} \;\; ({\rm the\;  Kronecker
\; delta}). $$ Now, statements $2-4$ follow.

$5$. For each $a\in A$, we have a finite sum $a=\sum \xba
\l_\alpha$ with $\l_\alpha \in A^\d$ (statement 1). Since $
\l_\alpha = \phi (\d^\alpha (a))$, we have $a= \sum\xba \phi
(\d^\alpha (a))$, and so statement 5.

$6$. The assumptions of statement $6$ guarantee that  $A_+$ is an
ideal of $A$, hence $\phi$ is an $A^\d$-algebra homomorphism.
$\Box$

{\bf Infinitely iterated Ore extensions}. Let $R$ be a ring, $d=\{
d_i, i\in I\}$ be a non-empty (not necessarily finite) set of
derivations of the ring $R$, $r=(r_{ij})$ be a skew symmetric
$I\times I$ matrix (possibly of infinite size) with entries from
$R$ ($r_{ij}=-r_{ji}$ and $r_{ii}=0$) such that
$$ [d_i, d_j] = \ad ( r_{ij})\;\; {\rm and} \;\; d_i(r_{jk})=-d_k(r_{ij})+d_j(r_{ik}) \;\; {\rm for \; all}\; i,j\in I. $$
For each {\em finite} subset, say $J=\{ 1, \ldots , n\}$, of $I$
consider the {\em iterated Ore extension}
$$ R_J:= R[t_1; d_1]\cdots [ t_n; d_n]= R_G[ t_n; d_n],\;\; G:= \{
1, \ldots , n-1\}, $$ where the derivation $d_n$ of the ring $R$
extended to a derivation (denoted in the same fashion) of the ring
$R_G$ by the rule $d_n(t_i)= r_{ni}$. It is an easy exercise  to
verify that $d_n$ is then well-defined: if $n=1$ then there is
nothing to prove. For $n>1$, by induction on $n$ one may assume
that $R_G$ is well-defined. Then the ring $R_G$ has the following
{\em defining} relations:
$$ [t_i, t_j]=r_{ij} \;\; {\rm and}\;\; [t_i, r]=d_i(r),  \;\; 0\leq
i<j\leq n-1, \; r\in R. $$ Now, the extended $d_n$ {\em respects}
these relations:
\begin{eqnarray*}
 d_n([t_i, t_j])&=& [d_n(t_i), t_j]+ [ t_i, d_n(t_j)]= [ r_{ni}, t_j] +[ t_i, r_{nj}]=-d_j (r_{ni})+d_i(r_{nj})= d_n(r_{ij}),  \\
 d_n([t_i, r])&=& [d_n(t_i), r]+ [ t_i, d_n(r)]= [ r_{ni}, r] +d_id_n(r)= ([d_n, d_i]+d_id_n)(r)= d_nd_i(r).
\end{eqnarray*}
Clearly, $J\subseteq L$ implies $R_J\subseteq R_L$ for finite
subsets $J, L\subseteq I$. The {\em infinitely iterated Ore
extension} $R[t ; d, r]$ is the union (the direct limit) of the
rings $R_J$. Without loss of generality one may assume that the
set $I$ is a well-ordered set, $(I, <)$. Then %\marginpar{Rts}
\begin{equation}\label{Rts}
R[t ; d, r]= \bigoplus_{\alpha \in \NI} Rt^\alpha
=\bigoplus_{\alpha \in \NI} t^\alpha R,
\end{equation}
where $t^\alpha$ is the {\em ordered } product
$t_{i_1}^{\alpha_1}\cdots t_{i_n}^{\alpha_n}$ where $\supp (\alpha
)=\{ \alpha_{i_1}, \ldots , \alpha_{i_n}\}$ and $ i_1<\ldots
<i_n$.

The {\em partial derivatives} $\der_i := \frac{\der}{\der t_i}\in
\Der_R(R[t ; d, r])$ are well-defined $R$-derivations of the ring
$R[t ; d, r]$ (since they respect the defining relations). So,
$\der := \{ \der_i, i\in I\}$ is the set of {\em commuting locally
nilpotent} $R$-derivations of the ring $R[t ; d, r]$ such that
$\der_i(t_j)=\d_{ij}$. If $R$ is an algebra over a field of
characteristic $p>0$ then $\der_i^p=0$, $i\in I$.

\begin{lemma}\label{f8Apr06}%\marginpar{f8Apr06}
Let the ring $\CA :=R[t ; d, r]$ be as above. Suppose, in
addition,  that $R$ is an algebra over a field of characteristic
$p>0$, $r_{ij}= [x_i, x_j]$ and $d_i=\ad (x_i)$ for some elements
$x_i\in R$. Then
\begin{enumerate}
\item $t_i^p-x_i^p$, $i\in I$, are central elements of $\CA$.
\item The ideal $\ga := (t_i^p-x_i^p, i\in I)$ of $\CA$ is
$\der_j$-invariant for all $j\in I$ ($\der_j (\ga )\subseteq
\ga$). \item Let $A:= \CA / \ga$ and $\d_i\in \Der_R(A)$: $a+\ga
\mapsto \der_i(a) +\ga$. Then $\d = \{ \d_i , i\in I\}$ is the set
of commuting $R$-derivations of $A$ such that $\d_i^p=0$, $\d_i
(x_j)=\d_{ij}$, and, for each $a\in A$, $\d_i(a)=0$ for almost all
$i$.
\end{enumerate}
\end{lemma}

{\it Proof}. $1$. $\ad (t_i^p)= \ad (t_i)^p= \ad (x_i)^p= \ad
(x_i^p)$, hence $\ad (t_i^p-x_i^p)=0$, i.e. $t_i^p-x_i^p$ is a
central element of $\CA$. The rest is obvious. $\Box $

The next theorem is the converse to Lemma \ref{f8Apr06}.

\begin{theorem}\label{8Apr06}%\marginpar{8Apr06}
Let $A$ be an algebra over a field $K$ of characteristic $p>0$,
$\d = \{ \d_i, i\in I\}$ be a non-empty set of commuting
$K$-derivations of the algebra $A$ such that $\d_i^p=0$,
$\d_i(x_j)=\d_{ij}$ (the Kronecker delta) for a set $x=\{ x_i,
i\in I\}$ of elements of $A$, and, for each $a\in A$, $\d_i (a)=0$
for almost all $i$. Then
\begin{enumerate}
\item $A= \bigoplus_{\alpha \in \CN } A^\d x^\alpha =
\bigoplus_{\alpha \in \CN }x^\alpha  A^\d $ where $\CN := \{
\alpha\in \NI \, | \,$ all $\alpha_i<p\}$ and $ x^\alpha :=
x_{i_1}^{\alpha_{i_1}}\cdots x_{i_n}^{\alpha_{i_n}}$ where $\supp
(\alpha ) =\{ \alpha_{i_1}, \ldots , \alpha_{i_n}\}$ and
$i_1<\cdots <i_n$ (where $(I,<)$ is the well-ordered set).  \item
$r_{ij}:=[x_i, x_j]\in A^\d$ and $[ x_i, A^\d ] \subseteq A^\d$
for all $i,j\in I$. \item $A\simeq A^\d [ t; d, r]/(t_i^p-x_i^p,
i\in I)$ where $t=\{ t_i, i\in I\}$, $d=\{ d_i:= \ad (x_i), i\in
I\}$, and $r=(r_{ij})$, $\{ t_i^p-x_i^p, i\in I\}$ are central
elements of the ring $A^\d [ t; d, r]$. \item Each derivation
$\d_i$ is induced by the partial $A^\d$-derivatives
$\der_i:=\frac{\der}{\der t_i}$ of the ring $A^\d [ t; d, r]$.
\end{enumerate}
\end{theorem}

{\it Proof}. $1$. Theorem \ref{9Apr06} since $E=\CN$.

$2$. For any $i,j\in I$, $ \d_j([x_i, A^\d ])\subseteq [ \d_{ij},
A^\d ] + [ x_i, \d_j(A^\d )]=0$, i.e. $[x_i, A^\d ] \subseteq
A^\d$. For $i,j,k\in I$, $\d_k ([x_i, x_j])=[ \d_{ki}, x_j]+[x_i,
\d_{kj}]=0$, i.e. $r_{ij}\in A^\d$.

$3$.  The ring $\CA :=A^\d [ t; d, r]$ is well-defined since
$[d_i, d_j] = [ \ad (x_i), \ad (x_j)]= \ad ([x_i, x_j])$ and $d_i
(r_{jk})= [x_i, [ x_j, x_k]]= [[x_i, x_j], x_k]+[x_j, [ x_i,
x_k]]= -d_k(r_{ij})+d_j(r_{ik})$, the {\em Jacobi} identity.

Each element $u:= t_i^p-x_i^p$ belongs to the centre of the ring
$\CA$ since it commutes with the generators of $\CA$: for each $\l
\in A^\d$, $[ t_i^p, \l ]= (\ad \, t_i)^p(\l)= (\ad \, x_i)^p(\l
)= [ x_i^p, \l ]$, and so $ [u, \l ] =0$; for each $j\in I$,
\begin{eqnarray*}
[ t_i^p, t_j] & =&(\ad \, t_i)^p(t_j)= (\ad \,
t_i)^{p-1}([t_i,t_j])=(\ad \, t_i)^{p-1}([x_i,x_j])=(\ad \,
x_i)^{p-1}([x_i,x_j])\\
&=&(\ad \, x_i)^p(x_j)=[x_i^p, x_j]=[x_i^p, t_j],\\
\end{eqnarray*}
and so $[u, t_j]=0$. $\Box $

%%%%%%%%%%%%%%%%%% SECTION 4 %%%%%%%%%%%%%%%%%%%%%%%%

\section{ A formula for the
inverse of an automorphism that preserves the ring of invariants}

Let $A$, $\d := \{ \d_i, i\in I\}$, and $ \{ \xbij \}$ be as in
Theorem \ref{9Apr06}.  Suppose that a $K$-algebra  automorphism
$\s\in \Aut_K(A)$  {\em preserves} the ring of invariants $A^\d$,
that is  $\s (A^\d )=A^\d$. Let $\s_\d := \s|_{A^\d}\in
\Aut_K(A^\d)$. Suppose that  we know explicitly the inverse
$\s_\d^{-1}$ and the twisted derivations
 $\d' := \{ \d_i':= \s \d_i\s^{-1}, i\in I\}$, then we can write down
 {\em explicitly} a formula for the
inverse automorphism $\s^{-1}$  (Theorem \ref{I9Apr06}). This is a
very simplified version of the  strategy we will follow in almost
all the examples later. Since $A=\oplus_{\alpha \in E}A^\d \xba =
\oplus_{\alpha \in E} \xba A^\d $ (Theorem \ref{9Apr06}), the
automorphism $\s$ is uniquely determined by its restriction
$\s_\d$ to the ring of invariants $A^\d$ and the images of the
iterative descents $ x_i':= \{ {x_i'}^{[j]}:= \s (\xbij )\}$,
$i\in I$. Note that the derivations $\d'=\{ \d_i'\}$ and their
iterative descents $ \{ x_i'\}$ do satisfy the conditions of
Theorem \ref{9Apr06} with $A^{\d'}= \s (A^\d )=A^\d$. For each
$i\in I$, let
\begin{eqnarray*}
\phi_i':=\sum_{k=0}^{l_i-1}(-1)^k\frac{(x_i')^k}{k!}(\d_i')^k, \;\;\;\;  \psi_i':=\sum_{k=0}^{l_i-1} (-1)^k(\d_i')^k (\cdot )\frac{(x_i')^k}{k!}, \\
\phi_\s := \phi':= \prod_{i\in I} \phi_i',  \;\;\;\;  \;\;\;\;
\;\;\;\;  \;\;\;\;  \;\;\;\;  \;\;\psi_\s := \psi':= \prod_{i\in
I} \psi_i',
\end{eqnarray*}
be the corresponding maps from Theorem \ref{9Apr06}. The maps $
\phi_\s , \psi_\s :A\ra A$ are {\em projections} onto the
subalgebra $A^{\d'}= A^\d$ of $A$, and,  for any $a\in A$,
%\marginpar{a4E}
\begin{equation}\label{a4E}
a=\sum_{\alpha \in E}{x'}^{[\alpha]} \phi_\s  ({\d'}^\alpha
(a))=\sum_{\alpha \in E}\psi_\s ({\d'}^\alpha (a)){x'}^{[\alpha
]}.
\end{equation}
 Then applying $\s^{-1}$ to these equalities we finish the
proof of the next theorem.

\begin{theorem}\label{I9Apr06}%\marginpar{I9Apr06}
Let $A$,  $\d =\{ \d_i\} $, $\d' =\{ \d_i'\} $, $\{ x_i\}$, $\{
x_i'\}$, and $\s $
 be as above. Then, for each $a\in A$,
$$ \s^{-1}(a)= \sum_{\alpha \in E} x^{[\alpha ]}\s_\d^{-1}(\phi_\s ({\d'}^\alpha
(a)))= \sum_{\alpha \in E}  \s_\d^{-1}(\psi_\s ({\d'}^\alpha
(a)))x^\alpha.$$
\end{theorem}

%%%%%%%%%%%%%%%%%% SECTION 5 %%%%%%%%%%%%%%%%%%%%%%%%

\section{An inversion formula for an automorphism of a central simple
algebra}\label{CCAIF}%\marginpar{CCAIF}

In this section,  $K$ is an {\em arbitrary} field and $A$ is a
{\em central simple} $K$-algebra (central means that the centre of
$A$ is $K$). For simplicity, let us assume that the algebra $A$ is
{\em countably generated}, for example $A$ is a finitely generated
algebra, the general case  is considered at the end of the
section. Then $A=\cup_{i\geq 0}A_i$ is a union of an ascending
chain of {\em finite dimensional}  subspaces: $K\subseteq
A_0\subseteq A_1\subseteq \cdots $.

The algebra $E:= {\rm End}_K(A)$ of all $K$-linear maps from $A$
to itself  is a {\em topological } algebra with respect to the
{\em finite topology} where neighbourhoods  of the zero map is
given by the ascending chain of subspaces
$$ E_0\supseteq E_1\supseteq \cdots \supseteq E_i:= \{ f\in E\, |
\, f(A_i)=0\} \supseteq \cdots , \;\; \cap_{i\geq 0} E_i=0.$$ So,
an element $f$ of $E$ is `small' if it annihilates `big' $A_i$.
Note that this description of the finite topology works only for
algebras that has  no more than countable basis over $K$. For the
general case, the reader is refereed to the book of  Jacobson
\cite{JacobsonSR}.

Let $\CA$ be the image of the $K$-algebra homomorphism
$$ A\t A^{op} \ra E, \;\; a\t b \mapsto (x\mapsto axb),$$
 where $x\in A$ and $A^{op}$ is the {\em opposite} algebra to $A$. By the {\em
 Density Theorem}, $\CA$ is a {\em dense} subalgebra in $E$, i.e.
 for any map $f\in E$ and any $i\geq 0$, there are elements $a_j,
 b_j\in A$ such that $f(a) =\sum_j a_j a b_j$ for all $a\in A_i$.

Let $e_1:=1, e_2, \ldots $ be a $K$-basis for $A$ such that for
each $i\geq 0$, $e_1=1, e_2, \ldots , e_{n_i}$ is a $K$-basis for
$A_i$. Clearly, $n_0\leq n_1\leq  \ldots $.  For each $i\geq 0$
and $j$ such that $1\leq j \leq n_i$, consider the $K$-linear map
$$ p_{ij}: A_i=\bigoplus_{k=1}^{n_i} Ke_k \ra A_i=\bigoplus_{k=1}^{n_i}
Ke_k, \;\; \sum_{k=1}^{n_i} \l_k e_k \mapsto \l_je_1= \l_j.$$
 By
the Density Theorem, there are elements $\{ a_{ij\nu}, b_{ij\nu},
\nu \in N_{ij}\}$ such that $p_{ij}(a)=\sum_{\nu \in
N_{ij}}a_{ij\nu } ab_{ij\nu}$ for all $a\in A_i$. Then, for each
$a\in A_i$, $a=\sum_{j=1}^{n_i} p_{ij} (a) e_j$, or, equivalently,
for each $a\in A$,
$$a=\sum_{j=1}^{n_i} p_{ij} (a) e_j, \;\; i\gg 0.$$
Applying $\s$ and denoting $\s (a)$ by $a$, we see that for each
$a\in A$,
$$a=\sum_{j=1}^{n_i} p_{ij}' (a) \s (e_j), \;\; i\gg 0,$$
where $p_{ij}' (\cdot )=\sum_{\nu \in N_{ij}}\s (a_{ij\nu}) (\cdot
) \s ( b_{ij\nu})$. Applying $\s^{-1}$, we prove the next theorem
(use the fact that $p_{ij}'(A)\subseteq K$).

\begin{theorem}\label{16Apr06}%\marginpar{16Apr06}
({\rm The inversion formula}) Let $A$ be the central simple
countably generated algebra over the field $K$ and $\s \in
\Aut_K(A)$. Then, for each $a\in A$,
$$ \s^{-1} (a) = \sum_{j=1}^{n_i} p_{ij}'(a) e_j, \;\; i\gg 0,$$
where $p_{ij}'(\cdot )= \sum_{\nu \in N_{ij}}\s (a_{ij\nu })
(\cdot ) \s (b_{ij\nu })$.
\end{theorem}

\begin{corollary}\label{c16Apr06}%\marginpar{c16Apr06}
({\rm The inversion formula for a central simple finite
dimensional algebra}) If, in addition, the algebra  $A$ is finite
dimensional over $K$. Then $A=A_i$ for some $i$ and,   for each
$a\in A$,
$$ \s^{-1} (a) = \sum_{j=1}^{n_i} p_{ij}'(a) e_j.$$
\end{corollary}

Let $A$ be an {\em arbitrary} central simple $K$-algebra. Fix a
covering $A=\cup_{i\in I} A_i$ of $A$ by a set of finite
dimensional subspaces $A_i$ such that $K\subseteq A_i$ for all
$i\in I$. For each $i\in I$, fix a basis $\{ e^i_j\}$ for the
vector space $A_i$ with $e^i_1:=1$, and then define the maps
$p_{ij}:A_i\ra A_i$  in the same way as before, $p_{ij}(\sum \l_k
e_k^i) =\l_je^i_1=\l_j\in K$. By the Density Theorem, there are
elements $\{ a_{ij\nu}, b_{ij\nu}\}$ such that $p_{ij}(a)=
\sum_{\nu} a_{ij\nu}ab_{ij\nu}$ for all $a\in A_i$.  Then each map
$p_{ij}':= \s (p_{ij})=\sum_{\nu} \s (a_{ij\nu}) (\cdot ) \s
(b_{ij\nu}):\s (A_i)\ra \s (A_i)$ has the image $K$ since $\s
(K)=K$, and, for each $a\in A$,
$$ a= \sum p_{ij}'(a) \s (e^i_j), \;\; i\gg 0,$$
where `$i\gg 0 $' means that for all $i$ such that $a\in \s
(A_i)$. Applying $\s^{-1}$, we have $\s^{-1} (a)= \sum p_{ij}'(a)
e^i_j$ for $i\gg 0$. This means that Theorem \ref{16Apr06} holds
with the new meaning of `$i\gg 0 $'.

%%%%%%%%%%%%%%%%%% SECTION 6 %%%%%%%%%%%%%%%%%%%%%%%%

\section{The inversion formula for an automorphism of a polynomial
algebra}\label{iforpol}%\marginpar{iforpol}

 In this section, $K$ is a field of characteristic $p>0$ and
 $P_n:=K[x_1, \ldots , x_n]$ is a polynomial algebra in $n$
 variables over the field $K$.

 Using the defining relations (\ref{DPndef}) for the algebra of differential operators  $\CD (P_n)$ on $P_n$ and its
 simplicity, one can verify that, {\em for each natural number $k\geq 1$, there exists a  $K$-algebra monomorphism of} $\CD (P_n)$ (which
is obviously {\bf not} {\em an isomorphism}) given by
 the rule
 %\marginpar{fkpk}
\begin{equation}\label{fkpk}
f_k: \CD (P_n)\ra \CD (P_n), \;\; x_i\mapsto x_i^{p^k}, \;\;\
 \der_i^{[\alpha ]}\ra \der_i^{[p^k\alpha ]}, \;\; i=1, \ldots ,
 n.
\end{equation}
 Note
that each monomorphism $f_k$, $k\geq 1$, is not a power
 of the {\em Frobenious map} $a\mapsto a^p$.

{\it Remark}. In the characteristic zero case (char$(K)=0$), it is
an old open question ({\it the Problem of Dixmier}, \cite{Dix},
Problem 1}): {\em is an algebra homomorphism of $\CD (P_n)$ an
algebra isomorphism?}

For an arbitrary field $K$, the ring $\CD (P_n)$ is {\em simple},
so any algebra homomorphism is automatically a {\em monomorphism}.
If ${\rm char} (K)=0$ then $\CD (P_n)=A_n$, but if ${\rm char}
(K)=p>0$ then $\CD (P_n)\neq A_n$ and the Weyl algebra $A_n$ has a
big centre and {\em because} of that the analogue of the Problem
of Dixmier for $A_n$  trivially {\em fails}:  if $A_1=K\langle x,
\der \, | \, \der x -x\der =1\rangle$ then the $K$-algebra
homomorphism $\alpha : A_1\ra A_1$, $ x\mapsto x+x^p$, $\der
\mapsto \der$, is not an automorphism since its restriction to the
centre $Z(A_1)= K[ x^p , \der ^p]$ (a polynomial algebra in two
variables) $ x^p\mapsto x^p+x^{p^2}$, $\der^p \mapsto \der^p$, is
not an automorphism. Though, the endomorphism $\alpha$ is a {\em
monomorphism}. We propose the following conjecture.

{\em Conjecture. Let $A_n$ be the Weyl algebra over a field $K$ of
characteristic $p>0$. Then each $K$-algebra homomorphism is a
monomorphism.}

The case $n=1$ is an easy consequence of the Tsen's Theorem
\cite{Tsuchi03}. For general $n$, I proved this conjecture for
certain {\em algebraic} extensions of the Weyl algebra $A_n$.

By (\ref{fkpk}), for each $i=1, \ldots , n$, and each $k\geq 0$,
$\der_i^{[p^k]}\in \Der_K(P_{n,i,k})$ (i.e. $
\der_i^{[p^k]}(P_{n,i,k})\subseteq P_{n,i,k}$),   $$ P_{n,i,k}:=
K[x_1, \ldots , x_{i-1}, x_i^{p^k}, x_{i+1}, \ldots , x_n], \;\;
P_{n,i,k}^{\der_i^{[p^k]}}=P_{n,i,k+1},\;\; (\der_i^{[p^k]})^p=0,
$$
and $\{ \frac{x_i^{p^kj}}{j!}, 0\leq j <p\}$ is the {\em iterative
$\der_i^{[p^k]}$-descent of maximal length} in $P_{n,i,k}$.

Hence, using step by step Corollary \ref{1pw2Dec06}, for each
$i=1, \ldots , n$, the map (an infinite product) %\marginpar{pipc1}
\begin{equation}\label{pipc1}
\phi_i:=\cdots \phi_{i,k}\cdots \phi_{i,1}\phi_{i,0}:P_n\ra P_n,
\;\; \phi_{i,k}:=
\sum_{j=0}^{p-1}(-1)^j\frac{x_i^{p^kj}}{j!}\der_i^{[p^kj]}, \;\;
k\geq 0,
\end{equation}
is a (well-defined)  homomorphism of $\Pnhi$-modules which is, in
fact, a {\em projection} onto the polynomial subalgebra $\Pnhi :=
K[x_1, \ldots , \widehat{x}_i, \ldots , x_n]$ (where $x_i$ is
missed) of $P_n$, i.e.
$$ \phi_i : P_n=\Pnhi \oplus P_nx_i\ra P_n=\Pnhi \oplus P_nx_i,
\;\; a+bx_i\mapsto a, \;\; {\rm where}\;\; a\in \Pnhi , \; b\in
P_n.$$ Note that, for a given polynomial $u\in P_n$ of degree  in
the variable $x_i$, say $d=d_0+d_1p+\cdots +d_kp^k$ where $0\leq
d_j<p$ and $d_k\neq 0$, one has $\phi_i(u)=\phi_{i,k}\cdots
\phi_{i,0}(u)$. For any fixed element $u\in P_n$, {\em almost all}
maps $\phi_{i,k}$ {\em in the product} (\ref{pipc1}) act as the
{\em identity map} on $u$. The maps $\phi_i$ {\em commute}, and
their product %\marginpar{pphis}
\begin{equation}\label{pphis}
\phi :=\prod_{i=1}^n\phi_i:P_n=K\oplus \sum_{i=1}^nP_nx_i\ra
P_n=K\oplus \sum_{i=1}^nP_nx_i, \;\; \sum_{\alpha \in
\Nn}\l_\alpha x^\alpha \mapsto \l_0, \;\; (\l_\alpha \in K)
\end{equation}
is a $K$-algebra homomorphism which is a {\em projection} onto the
field $K$. Let $u\in P_n$ and
$\deg_{x_i}(u)=d_{i,0}+d_{i,1}p+\cdots +d_{i,k_i}p^{k_i}$ be the
degree of the polynomial $u$ in $x_i$ where $0\leq d_{i,j}<p$.
Then  $\phi (u)=\prod_{i=1}^n \phi_{i,k_i}\cdots \phi_{i,0}(u)$,
i.e.  {\em almost all} maps $\phi_{i,k}$ in the product
(\ref{pphis}) act as the {\em identity} map on $u$.

Recall that, for $\alpha , \beta \in \mathbb{N}^n$,
%\marginpar{daaxb}
\begin{equation}\label{daaxb}
\dera (x^\beta ) = {\beta \choose \alpha } x^{\beta - \alpha},
\;\; {\beta \choose \alpha }:= \prod_i {\beta_i \choose \alpha_i},
\end{equation}
where, in the formula above, $x_i^t:=0$ for all {\em negative}
integers $t$ and all $i$.

The next result is a kind of the {\em Taylor} formula in {\em
prime characteristic}.

\begin{theorem}\label{11Jan06}%\marginpar{11Jan06}
For any $a \in P_n$,
$$ a=\sum_{\alpha \in \Nn}\phi (\derba (a))x^\alpha.$$
\end{theorem}

{\it Proof}. If $a=\sum \l_\alpha x^\alpha $, $\l_\alpha \in K$,
then, by (\ref{daaxb}) and (\ref{pphis}), $\phi (\derba
(a))=\l_\alpha$. $\Box $

Equivalently, the identity map ${\rm id}$ on $P_n$ can be written
as %\marginpar{pida}
\begin{equation}\label{pida}
{\rm id}(\cdot ) = \sum_{\alpha \in \Nn}\phi (\derba  (\cdot
))x^\alpha .
\end{equation}

Let $\s$ be  a $K$-automorphism  of the polynomial algebra $P_n$,
it is uniquely determined by the images of the canonical
generators, $x_1':= \s (x_1), \ldots , x_s':=\s (x_s)$. The {\em
Jacobian} of the automorphism $\s$, that is $\D := \det
(\frac{\der x_i'}{\der x_j})$, must be a {\em nonzero} scalar,
i.e. $\D  \in K^*:= K\backslash \{ 0\}$.  The  derivations $
\der_1':= \frac{\der }{\der x_1'}, \ldots , \der_s':= \frac{\der
}{\der x_s'} \in \Der_K(P_n)$ can be written as %\marginpar{pdad2}
\begin{equation}\label{pdad2}
\der_j'  := \D^{-1} \det
 \begin{pmatrix}
  \frac{\der \s (x_1)}{\der x_1} & \cdots & \frac{\der \s (x_1)}{\der x_n} \\
  \vdots & \vdots & \vdots \\
\frac{\der }{\der x_1} & \cdots & \frac{\der }{\der x_n}\\
 \vdots & \vdots & \vdots \\
\frac{\der \s (x_n)}{\der x_n} & \cdots & \frac{\der \s (x_n)}{\der x_n} \\
\end{pmatrix}, \;\;\; j=1, \ldots , n,
\end{equation}
where we have `dropped' $\s (x_j)$ in the determinant $\det
(\frac{\der \s (x_k)}{\der x_l})$.

 Let
$\dersba := \frac{{\der'}^\alpha }{\alpha !}:=
\frac{{\der_1'}^\alpha_1 }{\alpha !}\cdots \frac{{\der_n'}^\alpha
}{\alpha_n !}$  (where $\alpha =(\alpha_1, \ldots , \alpha_n)\in
\mathbb{N}^n$) be the corresponding {\em higher derivations} for
the {\em new} choice of  generators for the polynomial algebra
$P_n$, that is $x_1', \ldots , x_n'$. Now, we consider the
projections (\ref{pipc1}) and (\ref{pphis}) for the generators
$x_1', \ldots , x_n'$ of $P_n$. For each $i=1, \ldots , n$,
%\marginpar{1pipc1}
\begin{equation}\label{1pipc1}
\phi_i':=\cdots \phi_{i,k}'\cdots \phi_{i,1}'\phi_{i,0}':P_n\ra
P_n, \;\; \phi_{i,k}':=
\sum_{j=0}^{p-1}(-1)^j\frac{{x_i'}^{p^kj}}{j!}{\der_i'}^{[p^kj]},
\;\; k\geq 0.
\end{equation}
The maps $\phi_i'$ commute. Let %\marginpar{pdad4}
\begin{equation}\label{pdad4}
\phi_\s := \phi':= \prod_{i=1}^n \phi_i'.
\end{equation}

\begin{theorem}\label{i11Jan06}%\marginpar{i11Jan06}
{\rm (The Inversion Formula)} For each $\s \in \Aut_K(P_n)$
 and $a\in  P_n$,
 $$ \s^{-1}(a)=\sum_{\alpha \in \Nn}\phi_\s
 (\dersba (a))x^\alpha . $$
\end{theorem}

{\it Proof}. By Theorem \ref{11Jan06}, $a=\sum_{\alpha \in
\Nn}\phi_\s
 (\dersba (a)){x'}^\alpha$, then applying $\s^{-1}$
 we have the result
 $$ \s^{-1}(a)=\sum_{\alpha \in
\Nn }\phi_\s
 (\dersba (a))\s^{-1}({x'}^\alpha ) =\sum_{\alpha \in
\Nn }\phi_\s
 (\dersba (a))x^\alpha . \;\;\; \Box $$

%%%%%%%%%%%%%%%%%% SECTION 7 %%%%%%%%%%%%%%%%%%%%%%%%

\section{The inversion formula for  $\s \in
\Aut_K(\CD (K[x_1, \ldots ,
x_n]))$}\label{ifordop}%\marginpar{ifordop}

In this section, $K$ is a field of characteristic $p>0$ and $\CD
:= \CD (P_n)$ is the ring of differential operators on the
polynomial algebra $P_n:=K[x_1, \ldots , x_n]$ over the field $K$.
The algebra $\CD$ is a {\em central simple countably generated}
(but not finitely generated) algebra over $K$. The results of
Section \ref{CCAIF} show that the inversion formula for $\s \in
\Aut_K(\CD )$ can be written in the most economical way -  using
only addition and multiplication.  Clearly, $\CD =\CD (K[x_1])\t
\cdots \t \CD (K[x_n])$, the tensor product of algebras. For each
$i=1, \ldots , n$, the inner derivation $\d_i:= -\ad \, x_i$ of
the algebra $\CD $ is a {\em locally nilpotent} derivation,
$\CD^{\d_i} =\CD (K[x_1])\t\cdots \t K[x_i]\t  \cdots \t \CD
(K[x_n])$, and $\{ \derij , j\geq 0\}$ is the {\em iterative
$\d_i$-descent} (of {\em maximal} length) since
$$(-\ad \, x_i)(\der_i^{[j]})=\der_i^{[j-1]} \;\; {\rm and}\;\;
\der_i^{[j]}\der_i^{[k]}={j+k\choose j}\der_i^{[j+k]}\;\; {\rm for
\; all}\;\; j, k\geq 0.$$   By Theorem \ref{pw2Dec06}, there are
two projections onto $ \CD^{\d_i}$ (where $ a\in \CD^{\d_i}$):
\begin{eqnarray*}
 \phi_i=\sum_{k\geq 0}\derik (\ad \, x_i)^k:  \CD = \CD^{\d_i}\oplus \CD_{+,i}\ra  \CD^{\d_i}\oplus \CD_{+,i},
 a+b\mapsto a,
  b\in \CD_{+,i}= \oplus_{k\geq 1} \derik \CD^{\d_i}, \\
\psi_i=\sum_{k\geq 0} (\ad \, x_i)^k (\cdot ) \derik :  \CD =
\CD^{\d_i}\oplus \CD_{+,i}\ra \CD^{\d_i}\oplus \CD_{+,i},
 a+b\mapsto a,
  b\in \CD_{+,i}= \oplus_{k\geq 1}\CD^{\d_i} \derik .
\end{eqnarray*}
The maps $\phi_i$ and $\psi_i$ are homomorphisms of right and left
$ \CD^{\d_i}$-modules respectively. The maps $\phi_i$ (resp.
$\psi_i$) {\em commute} and their products yield projections onto
the polynomial subalgebra $P_n$ of the ring of differential
operators $\CD (P_n)$, %\marginpar{dph}
\begin{equation}\label{dph}
\phi_n\cdots \phi_1 :\CD  =P_n\oplus \CD_+\ra P_n\oplus \CD_+,\;
a+b\mapsto a, \; a\in P_n,\;
  b\in \CD_+:= \oplus_{0\neq \alpha \in \Nn} \derba P_n,
\end{equation}
%\marginpar{dps}
\begin{equation}\label{dps}
\psi_1\cdots \psi_n :\CD  =P_n\oplus \CD_+\ra P_n\oplus \CD_+,\;
a+b\mapsto a, \; a\in P_n,\;
  b\in \CD_+:= \oplus_{0\neq \alpha \in \Nn} P_n\derba .
\end{equation}

For each $i=1, \ldots , n$, and each $k\geq 0$, the
 inner derivation $\ad \, \der_i^{[p^k]}$ of $\CD $ preserves the subalgebra
$P_{n,i,k}:= K[x_1, \ldots , x_{i-1}, x_i^{p^k}, x_{i+1}, \ldots ,
x_n]$  (i.e. $ [\der_i^{[p^k]}, P_{n,i,k}] \subseteq P_{n,i,k}$),
$P_{n,i,k}^{\ad \, \der_i^{[p^k]}}=P_{n,i,k+1}$, $(\ad\,
\der_i^{[p^k]})^p=\ad (\der_i^{[p^k]})^p=\ad \, 0 = 0$, and $\{
\frac{x_i^{p^kj}}{j!}, 0\leq j <p\}$ is the {\em iterative $\ad \,
\der_i^{[p^k]}$-descent of maximal length} in $P_{n,i,k}$. By
Corollary \ref{1pw2Dec06}, the map

%\marginpar{Dpipc1}
\begin{equation}\label{Dpipc1}
\phi_{n+i}:=\cdots \phi_{n+i,k}\cdots
\phi_{n+i,1}\phi_{n+i,0}:P_n\ra P_n, \;\; \phi_{n+i,k}:=
\sum_{j=0}^{p-1}(-1)^j\frac{x_i^{p^kj}}{j!}(\ad\,
\der_i^{[p^k]})^j, \;\; k\geq 0,
\end{equation}
 is a
{\em projection} onto the subalgebra $\Pnhi :=K[x_1, \ldots ,
\widehat{x}_i, \ldots , x_n]$ of the polynomial algebra $P_n$,
$$ \phi_{n+i}:P_n=\Pnhi \oplus P_nx_i\ra P_n=\Pnhi \oplus P_nx_i,
\;\; a+bx_i\mapsto a, \;\; {\rm where}\;\; a\in \Pnhi , \; b\in
P_n.$$ The maps $\phi_{n+i}$ {\em commute} and their product
$$ \phi_{n+1}\cdots \phi_{2n}: P_n=K\oplus P_{n,+}\ra P_n=K\oplus
P_{n,+}, \;\; \sum_{\alpha \in \Nn} \l_\alpha x^\alpha \ra \l_0,
\;\; (\l_\alpha \in K),$$ is a $K$-algebra homomorphism which is a
projection onto the field $K$ where $P_{n,+}:= \sum_{i=1}^n
P_nx_i$. The maps $\phi_{n+i}$ are well-defined also as maps from
the algebra $\CD$ to itself as it  follows from the decomposition
$\CD =\oplus_{\alpha , \beta \in \Nn}Kx^\alpha \der^{[\beta ]}$
and (\ref{DPndef}). Note that the maps $\phi_i$, $\phi_j$ $(1\leq
i,j\leq 2n)$ commute unless $|i-j|=n$. The following maps (where
$a:=\sum\l_{\beta , \alpha} \derbb x^\alpha , a':=\sum\l_{ \alpha
, \beta }' x^\alpha\derbb \in \CD $, $\l_{\beta , \alpha} ,\l_{
\alpha , \beta }' \in K$): %\marginpar{CDphi1}
\begin{equation}\label{CDphi1}
\phi:=\phi_{2n}\cdots \phi_{n+1}\phi_n\cdots \phi_1 :\CD  =K\oplus
\CD_+\ra K\oplus \CD_+,\; a \mapsto \l_{0,0}, \; \CD_+:=
\bigoplus_{\alpha +\beta \neq 0 } K\derbb x^\alpha,
\end{equation}
%\marginpar{CDpsi1}
\begin{equation}\label{CDpsi1}
 \psi:=\phi_{2n}\cdots \phi_{n+1}\psi_1\cdots \psi_n :\CD  =K\oplus \CD_+\ra
K\oplus \CD_+,\; a' \mapsto \l_{0,0}', \; \CD_+:=
\bigoplus_{\alpha +\beta \neq 0 } Kx^\alpha \derbb ,
\end{equation}
are {\em projections} onto $K$.

The next result is a kind of a {\em noncommutative} Taylor formula
for the ring of differential operators $\CD (P_n)$ on $P_n$  in
{\em prime characteristic}.

\begin{theorem}\label{t12Jan06}%\marginpar{t12Jan06}
For any $a \in \CD (P_n)$,
$$ a=\sum_{\alpha , \beta \in \Nn}(-1)^{|\alpha |}\phi (\d^{ \beta} (a)\derba )\derbb x^\alpha =
\sum_{\alpha , \beta \in \Nn}\psi (\derba \d^{\beta}(
a))x^\alpha\derbb ,$$ where $\d^\beta :=\prod_{i=1}^n (-\ad \,
x_i)^{\beta_i}$ and $|\alpha | := \alpha_1+\cdots +\alpha_n$.
\end{theorem}

{\it Remark}. The element $\d^{ \beta} (a)\derba =\d^{ \beta} (a)
\cdot \derba $ (resp. $\derba \d^{\beta}( a)=\derba
\cdot\d^{\beta}( a)$) is the product in $\CD (P_n)$ of the
elements $\d^{ \beta}(a)$ and $ \derba$ (resp. $\derba $ and $
\d^{\beta}( a)$).

{\it Proof}. If $a=\sum \l_{\beta \alpha}\derbb  x^\alpha = \sum
\l_{ \alpha\beta}'  x^\alpha\derbb$ where  $\l_{\beta \alpha},
\l_{ \alpha\beta}' \in K$, then we must prove that $(-1)^{|\alpha
|}\phi (\d^{ \beta} (a)\derba )=\l_{\beta \alpha}$ and $\psi
(\derba \d^{\beta}( a))=\l_{\alpha \beta}'$. Since $\CD (P_n)=\CD
(K[x_1])\t\cdots \t\CD (K[x_n])$, it suffices to prove these
statements  when $n=1$ (for any $n$, repeat the arguments below
$n$ times or use induction on $n$). So, let $n=1$. Recall that
$\der_1^{[0]}:=1$ and $\der_1^{[s]}:=0$ for all {\em negative}
integers $s$.
\begin{eqnarray*}
 (-1)^i\phi (\d^j_1(a)\der_1^{[i]})&=& (-1)^i\phi ((-\ad \, x_1)^j
 (\sum_{\alpha, \beta \geq 0}\l_{\beta \alpha }\der_1^{[\beta ]}x_1^\alpha )\der_1^{[i]})
 = (-1)^i\phi (\sum_{\alpha \geq 0, \beta \geq j}\l_{\beta \alpha
 }\der_1^{[\beta-j ]}x_1^\alpha \der_1^{[i]}) \\
 &=&(-1)^i\phi (\sum_{\alpha \geq 0}\l_{j \alpha }x_1^\alpha
 \der_1^{[i]})=(-1)^i\phi (\sum_{\alpha \geq 0}\l_{j \alpha }\sum_{k=0}^\alpha {\alpha \choose k} (-1)^k
 \der_1^{[i-k]}x_1^{\alpha-k}) \\
 &=&(-1)^i(-1)^i\l_{ji}=\l_{ji},
\end{eqnarray*}
\begin{eqnarray*}
 \psi (\der_1^{[i]}\d^j_1(a))&=& \psi (\der_1^{[i]}(-\ad\,
 x_1)^j(a))= \psi  (\der_1^{[i]}(-\ad \, x_1)^j
 (\sum_{\alpha, \beta \geq 0}\l_{ \alpha \beta}'x_1^\alpha\der_1^{[\beta ]} ))\\
&=& \psi  (\der_1^{[i]}
 \sum_{\alpha \geq 0, \beta \geq j}\l_{ \alpha \beta}'x_1^\alpha\der_1^{[\beta -j]}
 )= \psi  (
 \sum_{\alpha \geq 0}\l_{ \alpha j}'\der_1^{[i]}x_1^\alpha
 )\\
 &=&\psi ( \sum_{\alpha \geq 0}\l_{ \alpha
j}' \sum_{k=0}^\alpha {\alpha \choose k}x_1^{\alpha -k}
\der_1^{[i-k]}) =\l_{ij}'. \;\; \Box
\end{eqnarray*}

Let  $\s $ be a $K$-automorphism of the algebra $\CD $. Then the
elements   $x_i':= \s (x_i)$, ${\der_i'}^{[k]}:= \s (\derik )$,
$i=1, \ldots , n, k\geq 1$, are another canonical generators for
the algebra $\CD$ (that satisfy the defining relations
(\ref{DPndef})).  Let $\phi_\s:= \phi'$ and $\psi_\s:= \psi'$ be
the maps as in  (\ref{CDphi1}) and (\ref{CDpsi1}) but for the new
canonical generators $x_i'$, ${\der_i'}^{[k]}$ (one has to put
$(')$ everywhere in the formulae).

\begin{theorem}\label{id12Jan06}%\marginpar{id12Jan06}
({\rm The
inversion formula for} $\s \in \Aut_K(\CD (P_n))$)
 For each $\s \in \Aut_K(\CD (P_n))$ and  $a \in
\CD (P_n)$,
$$ \s^{-1}(a)=\sum_{\alpha , \beta \in \Nn}(-1)^{|\alpha |}\phi_\s ({\d'}^{\beta} (a){\der'}^{|\alpha ]})\derbb x^\alpha =
\sum_{\alpha , \beta \in \Nn}\psi_\s ({\der'}^{ [\alpha ] }
 {\d'}^{\beta} (a))x^\alpha\derbb ,$$ where  ${\d'}^\beta
:=\prod_{i=1}^n (-\ad \, x_i')^{\beta_i}$.
\end{theorem}

{\it Proof}. By Theorem \ref{t12Jan06}, $$a=\sum_{\alpha , \beta
\in \Nn}(-1)^{|\alpha |}\phi_\s ({\d'}^{\beta} (a){\der'}^{|\alpha
]}){\der'}^{[\beta ]} {x'}^\alpha = \sum_{\alpha , \beta \in
\Nn}\psi_\s ({\der'}^{ [\alpha ] }
 {\d'}^{\beta} (a)){x'}^\alpha {\der'}^{[\beta
]}.$$ Now, applying $\s^{-1}$ to these equalities  we have the
result. $\Box $

%%%%%%%%%%%%%%%%%% SECTION 8 %%%%%%%%%%%%%%%%%%%%%%%%

\section{The inversion formula for  $\s \in
\Aut_K(\CD (P_n)\t P_m)$}\label{ifornm}%\marginpar{ifordnm}

In this section, $K$ is a field of characteristic $p>0$ and $\CD
:= \CD (P_n)$ is the ring of differential operators on the
polynomial algebra $P_n:=K[x_1, \ldots , x_n]$ over the field $K$,
$P_m:= K[y_1, \ldots , y_m]$ is a polynomial algebra in $m$
variables over $K$, and $A:= \CD (P_n)\t P_m$. The inversion
formula (Theorem \ref{i13Jan06}) for $\s \in \Aut_K(A)$ is a
consequence of the results of the previous two sections, very
briefly, the main reason for that is that the algebra $\CD (P_n)$
is {\em central}, i.e. the centre of the algebra $\CD (P_n)$ is
$K$. Therefore, the centre of the algebra $A$ is $P_m$, and so $\s
(P_m)=P_m$. Note that $ A=\bigoplus_{\alpha , \beta \in \Nn ,
\g\in \Nm} Kx^\alpha \derbb y^\g = \bigoplus_{\alpha , \beta \in
\Nn , \g\in \Nm} K\derbb x^\alpha  y^\g $.

The maps $\phi $ and  $\psi$ from (\ref{CDphi1}) and
(\ref{CDpsi1}) respectively  make sense for the algebra $A$ by
simply extending `scalars' in the natural way ($ \CD_{Q_m}(P_m\t
Q_m)\simeq \CD_K(P_n)\t Q_m$ where $Q_m:= K(y_1, \ldots , y_m)$ is
the field of rational functions). Let us denote them by $\phi_\CD
$ and  $\psi_\CD $ respectively.  The maps
\begin{eqnarray*}
\phi_\CD :& A=P_m\oplus A_+\ra P_m\oplus A_+, \; a\mapsto \l_0, \;
A_+:=\bigoplus_{|\alpha |+ |\beta | \neq
0} P_m\derbb x^\alpha  , \\
\psi_\CD :& A=P_m\oplus A_+\ra P_m\oplus A_+, \; a'\mapsto \l_0',
\; A_+:=\bigoplus_{|\alpha |+ |\beta | \neq 0} P_mx^\alpha \derbb
 ,
\end{eqnarray*}
are projections onto the centre  $P_m$ of $A$ where $a=\sum \l_{
\beta  \alpha } \derbb x^\alpha $, $a'=\sum \l_{ \alpha \beta }'
x^\alpha \derbb  $, $\l_{ \beta  \alpha }, \l_{ \alpha  \beta
}'\in P_m$. Let $\phi_{P_m}: P_m\ra P_m$ be the map (\ref{pphis}).

Now, the next theorem is a direct consequence of Theorem
\ref{11Jan06} and Theorem \ref{t12Jan06}.

\begin{theorem}\label{t13Jan06}%\marginpar{t13Jan06}
For any $a \in A:=\CD (P_n)\t P_m$,
\begin{eqnarray*}
 a&=&\sum_{\alpha , \beta \in
\Nn , \g \in \Nm}(-1)^{|\alpha |}\phi_{P_m} (\der^{[\g ]}
(\phi_\CD (\d^\beta (a)\dera ))) \derbb x^\alpha y^\g \\
 &=&\sum_{\alpha , \beta \in \Nn , \g \in \Nm}\phi_{P_m} (\der^{[\g ]}
(\psi_\CD (\dera \d^\beta(a))))x^\alpha\derbb y^\g ,
\end{eqnarray*}
 where
$\d^{\beta}$ is as in Theorem \ref{t12Jan06}, and $\der^{[\g ]}$
is  as in Theorem \ref{11Jan06}.
\end{theorem}

Let  $\s \in \Aut_K(A)$. Then its restriction to the centre
$Z(A)=P_m$ of the algebra $A$ is an automorphism, say  $\tau \in
\Aut_K(P_m)$. Let $\phi_{P_m, \s}$ be the map (\ref{pdad4}) for
the restriction $\tau$. Similarly, let $\phi_{\CD , \s}$ and
$\psi_{\CD , \s}$ be the maps given by the same formulae as in
 Theorem \ref{id12Jan06} but for the algebra $A=\CD \t P_m$ rather
than $\CD $, let us stress it again that  this means that we
extend the `scalars' from $K$ to $P_m$ (of course, one can extend
the field from $K$ to the field of rational functions $Q_m:=
K(x_1, \ldots , x_m)$ and repeat all the arguments from Section
\ref{ifordop} over this bigger field). One should stress that the
formulae for $\phi_{\CD , \s}$ and $\psi_{\CD , \s}$ are the same,
the only new thing is that `scalars'  of all the elements and of
the inner derivations involved are in $P_m$ rather than $K$.

\begin{theorem}\label{i13Jan06}%\marginpar{i13Jan06}
({\rm The
inversion formula for} $\s \in \Aut_K(\CD (P_n)\t P_m)$)
 For each $\s \in \Aut_K(\CD (P_n)\t P_m)$ and  $a \in
\CD (P_n)\t P_m$,
\begin{eqnarray*}
 \s^{-1}(a)&=& \sum_{\alpha , \beta \in
\Nn , \g \in \Nm}(-1)^{|\alpha |}\phi_{P_m, \s } ({\der'}^{[\g ]}
(\phi_{\CD , \s } ({\d'}^\beta (a){\der'}^{[\alpha ]} ))) \derbb x^\alpha y^\g \\
 &=&\sum_{\alpha , \beta \in \Nn , \g \in \Nm}\phi_{P_m, \s } ({\der'}^{[\g ]}
(\psi_{\CD , \s } ({\der'}^{[\alpha ]}
{\d'}^\beta(a))))x^\alpha\derbb y^\g .
\end{eqnarray*}

\end{theorem}

{\it Proof}. By Theorem \ref{t13Jan06},
\begin{eqnarray*}
a&=& \sum_{\alpha , \beta \in \Nn , \g \in \Nm}(-1)^{|\alpha
|}\phi_{P_m, \s } ({\der'}^{[\g ]}
(\phi_{\CD , \s } ({\d'}^\beta (a){\der'}^{[\alpha ]} ))) {\der'}^{[\beta ]} {x'}^\alpha {y'}^\g \\
 &=&\sum_{\alpha , \beta \in \Nn , \g \in \Nm}\phi_{P_m, \s } ({\der'}^{[\g ]}
(\psi_{\CD , \s } ({\der'}^{[\alpha ]}
{\d'}^\beta(a)))){x'}^\alpha{\der'}^{[\beta ]} {y'}^\g .
\end{eqnarray*}
Applying $\s^{-1}$ yields the result.  $\Box $

{\it Remark}. It makes sense to stress that basically the process
of finding the  inversion formula for an automorphism of the
algebra $A:=\CD \t P_m$ collapses to two cases - the polynomial
and the ring of differential operator's case - simply because the
ring of differential operators $\CD (P_n)$ is a {\em central}
algebra. We will see later that for the Weyl algebra $A_n$ this is
not the case - the centre of $A_n$ is big, it is a polynomial
algebra in $2n$ variables, and the problem of finding the
inversion formula for $\s \in \Aut_K(A_n\t P_m)$ can not be
reduced to the cases : $\s \in \Aut_K(A_n)$ and $\s \in
\Aut_K(P_m)$ (see Section \ref{iforWeyl}) in such a
straightforward manner as in this case.

%%%%%%%%%%%%%%%%%% SECTION 9 %%%%%%%%%%%%%%%%%%%%%%%%

\section{The inversion formula for an automorphism of the $n$'th Weyl algebra $A_n$ (and of $A_n\t K[y_1, \ldots ,
y_m]$)}\label{iforWeyl}%\marginpar{iforWeyl}

Let $K$ be a field of characteristic $p>0$. The $n$'th {\em Weyl }
algebra $A_n=A_n(K)$ is a $K$-algebra generated by $2n$ generators
$q_1, \ldots , q_n$, $p_1, \ldots , p_n$ which are subject to the
defining relations:
$$ [p_i, q_j]=\d_{ij}, \;\; [p_i, p_j]=[q_i, q_j]=0\;\; {\rm
for\;\; all}\;\; i,j=1, \ldots , n,$$ where $\d_{ij}$ is the
Kronecker delta, $[a,b]:=ab-ba$.
 The Weyl algebra $A_n=\oplus_{\alpha , \beta \in \Nn} Kp^\beta
 q^\alpha =\oplus_{\alpha , \beta \in \Nn} K
 q^\alpha p^\beta $ (where $p^\beta := p_1^{\beta_1}\cdots
 p_n^{\beta_n}$ and $q^\alpha := q_1^{\alpha_1}\cdots
 q_n^{\alpha_n}$) is a Noetherian algebra which is a free finitely
 generated module over its centre $Z(A_n)=K[x_1, \ldots , x_{2n}]$,
 the
 polynomial algebra in $2n$ variables
 $$ x_1:=q_1^p, \ldots , x_n:=q_n^p, x_{n+1}:=p_1^p, \ldots ,
 x_{2n}:= p_n^p.$$
 Clearly, $A_n=\oplus_{\alpha , \beta \in \CN}Z(A_n)p^\beta
 q^\alpha =\oplus_{\alpha , \beta \in \CN}Z(A_n)q^\alpha p^\beta$
 where $\CN := \{ \alpha \in \Nn \, | \, 0\leq \alpha_1<p, \ldots
 , 0\leq \alpha_n<p\}$, and $A_n=K\langle p_1, q_1\rangle \t
 \cdots \t K\langle p_n, q_n\rangle \simeq A_1^{\t n }$, the
 tensor product of $n$ copies of  the first Weyl algebra. Let $P_m:= K[x_{2n+1},
 \ldots , x_{2n+m}]$ be a polynomial algebra in $m$ variables and
 $P_0:=K$. Let $A:= A_n\t P_m$ be the tensor product of algebras.
 In particular, $A=A_n$ if $m=0$. The generators $q_1, \ldots ,
 q_n, p_1, \ldots ,
 p_n, x_{2n+1}, \ldots , x_{2n+m}$ will be called the {\em
 canonical } generators of the algebra $A$. Let $r$ be one of the
 canonical generators of the {\em Weyl  algebra} $A_n$ (i.e. $r\neq x_{2n+i}$ for all $i$)
  and let   $\Ahr
 $ be the subalgebra of $A$ generated by $r^p$ and all its canonical
 generators except $r$. Clearly, $ (\ad \, r)^p=\ad (r^p)=0$
 since the element $r^p$ belongs to the centre $Z$ of the algebra
 $A$ (where $\ad \, r$ is the inner derivation of the algebra
 $A$). By Corollary \ref{1pw2Dec06}, for each $i=1, \ldots  , n$,  the
 maps
 \begin{eqnarray*}
\Phi_i := & \sum_{j=0}^{p-1}(-1)^j\frac{q_i^j}{j!}(\ad \, p_i)^j:
A=\oplus_{j=0}^{p-1}q_i^j A_{\widehat{q}_i}\ra A, \\
\Psi_i := & \sum_{j=0}^{p-1}(-1)^j(\ad \, p_i)^j(\cdot )
\frac{q_i^j}{j!}: A=\oplus_{j=0}^{p-1}A_{\widehat{q}_i}q_i^j\ra A,
 \end{eqnarray*}
are projections onto the subalgebra $A_{\widehat{q}_i}$ of $A$,
and the maps
 \begin{eqnarray*}
\Phi_{n+i} := & \sum_{j=0}^{p-1}\frac{p_i^j}{j!}(\ad \, q_i)^j:
A=\oplus_{j=0}^{p-1}p_i^j A_{\widehat{p}_i}\ra A, \\
\Psi_{n+i} := & \sum_{j=0}^{p-1}(\ad \, q_i)^j(\cdot )
\frac{p_i^j}{j!}: A=\oplus_{j=0}^{p-1}A_{\widehat{p}_i}p_i^j\ra A,
 \end{eqnarray*}
are projections onto the subalgebra $A_{\widehat{p}_i}$ of $A$.
Then their compositions %\marginpar{1Wphil}
\begin{equation}\label{1Wphil}
\Phi := \Phi_{2n}\Phi_{2n-1}\cdots \Phi_1: A=\oplus_{\alpha ,
\beta \in \CN}Zq^\alpha p^\beta \ra A,
\end{equation}
%\marginpar{1Wpsil}
\begin{equation}\label{1Wpsil}
\Psi := \Psi_{2n}\Psi_{2n-1}\cdots \Psi_1: A=\oplus_{\alpha ,
\beta \in \CN}Zp^\beta q^\alpha \ra A,
\end{equation}
are projections onto the centre $Z:=Z(A_n)\t P_m$ of the algebra
$A$. The centre $Z$ is a polynomial algebra $P_s=K[x_1, \ldots ,
x_{2n}, x_{2n+1},\ldots x_s]$ in $s:= 2n+m$ variables, let
$\phi_Z$ be the map as in  (\ref{pphis}) in the case of the
polynomial algebra $Z=P_s$. Then the maps %\marginpar{Wphi1}
\begin{equation}\label{Wphi1}
\phi :=\phi_Z\Phi : A\ra A, \;\; a=\sum_{\alpha , \beta \in \CN ,
\g \in \Ns} \l_{\alpha \beta \g } q^\alpha p^\beta x^\g \mapsto
\l_0, \;\; (\l_{\alpha \beta \g }\in K),
\end{equation}
%\marginpar{Wpsi1}
\begin{equation}\label{Wpsi1}
\psi :=\phi_Z\Psi : A\ra A, \;\; a=\sum_{ \alpha , \beta  \in \CN
, \g \in \Ns} \l_{ \beta \alpha  \g }  p^\beta q^\alpha x^\g
\mapsto \l_0, \;\; (\l_{ \beta\alpha \g }\in K),
\end{equation}
are projections onto the field $K$.  The inner derivations $\ad\,
q_1, \ldots , \ad\, q_n, \ad\, p_1,\ldots , \ad \, p_n$ of the
algebra $A$ {\em commute}.
\begin{theorem}\label{tW14Jan06}%\marginpar{tW14Jan06}
For any $a \in A$,
$$ a=\sum_{\alpha , \beta \in \CN , \g \in \Ns}\phi_Z (\der^{[\g ]} (\Phi ( \d^{[\alpha , \beta ]}  (a)))) q^\alpha p^\beta x^\g  =
\sum_{\alpha , \beta \in \CN , \g \in \Ns}\phi_Z (\der^{[\g ]}
(\Psi (\d^{[\alpha , \beta ]}  (a))) p^\beta q^\alpha  x^\g ,$$
where
 $\der^{[\g ]}:= \prod_{i=1}^s\der_i^{[\g_i]}$,  $\d^{[\alpha , \beta ]}:=
 (\alpha !\beta !)^{-1}\prod_{i=1}^n (\ad \, p_i)^{\alpha_i}\prod_{j=1}^n (-\ad \, q_j)^{\beta_j}$.
\end{theorem}

{\it Proof}. If $a= \sum \l _{\alpha \beta \g }q^\alpha p^\beta
x^\g =\sum \l _{ \beta\alpha \g }'p^\beta q^\alpha  x^\g\in A$
where $\l _{\alpha \beta \g },\l _{ \beta\alpha \g }'\in K$, then
$\phi_Z (\der^{[\g ]}(\Phi ( \d^{[\alpha , \beta ]}
(a)))=\l_{\alpha \beta \g}$ and $\phi_Z (\der^{[\g ]} (\Psi
(\d^{[\alpha , \beta ]}  (a)))=\l _{ \beta\alpha \g }'$. $\Box$

Let $\s\in \Aut_K(A)$ and let $q_i':=\s (q_i)$, $p_i':= \s(p_i)$,
$i=1,\ldots , n$, and $x_j':= \s (x_j)$, $j=1, \ldots , s$. Then
the elements $q_1', \ldots ,q_n', p_1', \ldots ,   p_n',
x_{2n+1}', \ldots , x_s'$ are another choice of the canonical
generators for the algebra $A$ and $x_1', \ldots , x_s'$ are
generators for the polynomial algebra $Z=P_s$. Let $\phi_{Z, \s
}$, $\Phi_\s $ and $\Psi_\s$ be the maps (\ref{pphis}),
(\ref{1Wphil}) and (\ref{1Wpsil}) for the new choice of the
canonical generators: for, we have to put $(')$ everywhere.

\begin{theorem}\label{iW14Jan06}%\marginpar{iW14Jan06}
({\rm The
inversion formula for} $\s \in \Aut_K(A)$)  For any $\s \in
\Aut_K(A)$ and  $a \in A$,
\begin{eqnarray*}
\s^{-1} (a)& = & \sum_{\alpha , \beta \in \CN , \g \in
\Ns}\phi_{Z, \s} ({\der'}^{[\g ]} (\Phi_\s ( {\d'}^{[\alpha ,
\beta ]}  (a))))
 q^\alpha p^\beta x^\g  \\
 &=& \sum_{\alpha , \beta \in \CN , \g \in \Ns}\phi_{Z, \s}
({\der'}^{[\g ]} (\Psi_\s ( {\d'}^{[\alpha , \beta ]}  (a))))
p^\beta q^\alpha x^\g ,
\end{eqnarray*}
 where
 ${\der'}^{[\g ]}:= \prod_{i=1}^{2n}{\der_i'}^{[\g_i]}$,  ${\d'}^{[\alpha , \beta ]}:=
 (\alpha !\beta !)^{-1}\prod_{i=1}^n (\ad \, p_i')^{\alpha_i}\prod_{j=1}^n (-\ad \, q_j')^{\beta_j}$.
\end{theorem}

{\it Proof}. By Theorem \ref{tW14Jan06},
\begin{eqnarray*}
 a&=& \sum_{\alpha , \beta \in \CN , \g \in \Ns}\phi_{Z, \s}
({\der'}^{[\g ]} (\Phi_\s ({\d'}^{[\alpha , \beta ]}  (a))))
{q'}^\alpha {p'}^\beta {x'}^\g \\
 &=&  \sum_{\alpha , \beta \in \CN , \g
\in \Ns}\phi_{Z, \s} ({\der'}^{[\g ]} (\Psi_\s ({\d'}^{[\alpha ,
\beta ]} (a)))) {p'}^\beta {q'}^\alpha  {x'}^\g ,
\end{eqnarray*}
 then applying
$\s^{-1}$ proves the result.  $\Box$

%%%%%%%%%%%%%%%%%% SECTION 10 %%%%%%%%%%%%%%%%%%%%%%%%

\section{The inversion formula for  $\s \in \Aut_{K,c}(K[[x_1, \ldots , x_n]]$}

The notations of Section \ref{iforpol} remain fixed. Let $\hP_n$
be a series algebra in $n$ variables $x_1, \ldots , x_n$ over a
field $K$ of characteristic $p>0$. The algebra $\hP_n$ is a
completion of the polynomial algebra $P_n:=K[x_1, \ldots , x_n]$
with respect to the $\gm $-adic topology given by the powers of
the maximal ideal $\gm := (x_1, \ldots , x_n)$ of the algebra
$P_n$. So, the algebra $\hP_n$ is a complete local algebra with
maximal ideal $\hgm := \hP_n \gm$.

The higher derivations $\derij\in \CD (P_n)$ of the polynomial
algebra $P_n$ are {\em continuous} maps in the $\gm$-adic
topology, $\derij (\gm^k)\subseteq \gm^{k-j}$ for all $k\geq 0$
where $\gm^{-l}:=P_n$ for $l\geq 0$. Hence, so are the maps
$\phi_{i,k}$ from (\ref{pipc1}). It follows from the definition of
the maps $\phi_i$ (see (\ref{pipc1})) that they are well-defined
and continuous in the $\gm$-adic topology, hence so is the map
$\phi$ from (\ref{pphis}) as their finite product. We denote by
the same symbols $\phi_{i,k}$, $\phi_i$, $\phi$, $\derij$,
$\derba$, etc, the unique extensions of these continuous maps to
the (necessarily continuous) maps from $\hP_n$ to itself (note
that $\phi_{i,k}\cdots \phi_{i,0}(\sum_{j=0}^{p^{k+1}-1}
a_jx_i^j)=a_0$ where $a_j\in K[[ x_1, \ldots , x_{i-1} ,
x_i^{p^{k+1}}, x_{i+1}, \ldots
 , x_n]]$). The maps $\phi_i$ commute, and the map %\marginpar{spphis}
\begin{equation}\label{spphis}
\phi : \hP_n=K\oplus \hgm \ra K\oplus \hgm , \;\; \sum_{\alpha \in
\Nn}\l_\alpha x^\alpha \mapsto \l_0, \;\;(\l_\alpha \in K)
\end{equation}
is \ projection onto $K$.
\begin{theorem}\label{t14Jan06}%\marginpar{t14Jan06}
For any $a \in \hP_n$,
$$ a=\sum_{\alpha\in \Nn}\phi (\derba   (a)) x^\alpha .$$
\end{theorem}

{\it Proof}. If $a= \sum \l _{\alpha } x^\alpha $, $\l_\alpha \in
K$, then $\phi (\derba (a))=\l_\alpha $. $\Box$

Let $\s :\hP_n\ra \hP_n$ be a continuous $K$-algebra {\em
endomorphism} such that $\s (\hgm )\subseteq \hgm $ (this is, in
fact, a part of the definition of  continuous endomorphism) and
such that  its {\em Jacobian} $\D := \det (\frac{\der x_i'}{\der
x_j})$ is a {\em unit} of the algebra $\hP_n$ where $x_1':= \s
(x_1),\ldots , x_n':= \s (x_n)$. Then obviously $\s \in
\Aut_{K,c}(\hP_n)$ where $\Aut_{K,c}(\hP_n)$ is the group of {\em
continuous automorphisms} of the algebra $\hP_n$ (if $\tau \in
\Aut_{K,c}(\hP_n)$ then $\tau (\gm )\subseteq \gm$, by
definition).

{\it Example}. Let $\s$ be a $K$-algebra {\em endomorphism} of the
polynomial algebra $P_n$ such that $\s (\gm ) \subseteq \gm$ and
its Jacobian $\D \in K^*$. Then the $\s$ can be extended uniquely
to a continuous $K$-automorphism of the algebra $\hP_n$.

Given a {\em continuous} automorphism of the algebra $\hP_n$ then
its Jacobian is automatically  a unit of the algebra $\hP_n$ as
follows immediately from the chain rule.

So, let $\s \in \Aut_{K, c}(\hP_n)$. Let us define  continuous
maps $\der_j'$, ${\der'}^{[\alpha ]}$, $\phi_i'$, $\phi_\s
:\hP_n\ra \hP_n$ in the same way as in (\ref{pdad2}),
(\ref{1pipc1}), and (\ref{pdad4}) respectively.

\begin{theorem}\label{i14Jan06}%\marginpar{i14Jan06}
({\rm The
inversion formula for} $\s \in \Aut_{K,c}(\hP_n)$) For any $ \s
\in \Aut_{K,c}(\hP_n)$ and  $a \in \hP_n$,
$$ \s^{-1}(a)=\sum_{\alpha\in \Nn}\phi_\s ({\der'}^{[\alpha ]}   (a)) x^\alpha .$$
\end{theorem}

{\it Proof}. By Theorem \ref{t14Jan06},  $a=\sum_{\alpha\in
\Nn}\phi_\s ({\der'}^{[\alpha ]}   (a)) {x'}^\alpha$, then
applying $\s^{-1}$ we have the result. $\Box$

%%%%%%%%%%%%%%%%%% SECTION 11 %%%%%%%%%%%%%%%%%%%%%%%%

\section{The inversion formula for $\s \in \Aut_K(T_{k_1,
\ldots , k_n}\t P_m)$}\label{IFTk}%\marginpar{IFTk}

Let $K$ be a field of characteristic $p>0$, $\CD := \CD (P_n)$ be
the ring of differential operators on the polynomial algebra
$P_n:= K[ x_1, \ldots , x_n]$. For each $\k =(k_1, \ldots ,
k_n)\in \mathbb{N}^n$, consider the $K$-subalgebra $T_\k :=
T_{k_1, \ldots , k_n}$ of $\CD $ generated by the polynomial
algebra $P_n$ and the elements $\der_i^{[p^{j_i}]}$, $i=1, \ldots
, n$, $j_1<k_1, \ldots , j_n<k_n$.  The algebras $T_\k$ play an
important role in studying the ring $\CD $ and its modules. Then
$$ T_\k =\bigoplus_{\alpha \in \CN} P_n \dera =\bigoplus_{\alpha \in \CN}  \dera P_n,$$
where $\CN := \{ \alpha = (\alpha_i)\in \mathbb{N}^n\, | \,
\alpha_1<p^{k_1}, \ldots , \alpha_n<p^{k_n}\}$. One can easily
verify that the algebra $T_{1, \ldots , 1}$ is canonically
isomorphic to the factor algebra $A_n/(p_1^p, \ldots , p_n^p)$ of
the Weyl algebra $A_n$.  Note that if $k_i=0$ then there is no
element $\der_i^{[p^{j}]}$ with $j<0$. In order to accommodate
this border case, we will assume that $k_1\geq 1, \ldots , k_n\geq
1$  but instead of the algebra $T_\k$ we consider the algebra
$$ T:= T_\k \t P_m, \;\; P_m:= K[x_{n+1}, \ldots x_{n+m}], \;\;
s:= n+m,$$ (it is obvious that for a general $\k$, the algebra
$T_\k$ is of this type). Let $P_s:= K[x_1, \ldots , x_s]$.  Then
the algebra
$$ T= \bigoplus_{\alpha \in \CN} P_s\dera = \bigoplus_{\alpha \in \CN}
\dera P_s$$ is a left and right free finitely generated
$P_s$-module of rank $p^{k_1+\cdots +k_n}$ with the centre $Z:=
Z(T)= K[x_1^{p^{k_1}}, \ldots , x_n^{p^{k_n}}, x_{n+1}, \ldots ,
x_s]$, a polynomial algebra in $s$ indeterminates. The algebra
%\marginpar{TZab}
\begin{equation}\label{TZab}
T= \bigoplus_{\alpha, \beta \in \CN} Zx^\alpha \derb =
\bigoplus_{\alpha, \beta \in \CN} Z \derb x^\alpha
\end{equation}
is a  free finitely generated $Z$-module of rank $p^{2(k_1+\cdots
+k_n)}$. Let us consider the set $\d := \{ \d_i := -\ad \, x_i,
i=1, \ldots , n\}$ of {\em commuting nilpotent} $K$-derivations of
the algebra $T$, $\d_i^{p^{k_i}}= (-\ad \, x_i)^{p^{k_i}}=
(-1)^{p^{k_i}}\ad \, x_i^{p^{k_i}}=0$ since $x_i^{p^{k_i}}\in Z$.
The set $\der_i^*:= \{ \derij , 0\leq j <p^{k_i}\}$ is the {\em
iterative $\d_i$-descent} of maximal length. Clearly, $T^\d= P_s$.
The set $\d$ satisfies the conditions of Theorem \ref{9Apr06}, let
\begin{eqnarray*}
 \prod_{i=1}^n \phi_i : T\ra T, & \phi_i:= \sum_{k=0}^{p^{k_i}-1} \derik (\ad \, x_i)^k,  \\
  \prod_{i=1}^n \psi_i  : T\ra T, & \psi_i:= \sum_{k=0}^{p^{k_i}-1}  (\ad \, x_i)^k(\cdot )\derik,
\end{eqnarray*}
be the corresponding maps from (\ref{psphb}), these are the
projections onto $P_s$ as in Theorem \ref{9Apr06}.(3).

For each $i=1, \ldots , n$ and each $k=0, \ldots , k_i-1$, the
inner derivation $\ad \, \der_i^{[p^k]}$ of the algebra $T$ is a
nilpotent derivations: $(\ad \, \der_i^{[p^k]})^p=\ad \,
(\der_i^{[p^k]})^p=\ad \, 0 =0$, and the subalgebra of $P_s$,
$$ P_{s,i,k}:= K[ x_1, \ldots , x_{i-1}, x_i^{p^k}, x_{i+1},
\ldots , x_n]\t P_m$$ is $\ad \, \der_i^{[p^k]}$-{\em invariant},
the kernel of the derivation $\ad \, \der_i^{[p^k]}$ in
$P_{s,i,k}$ is equal to $P_{s,i,k+1}$. Note that $\{
\frac{x_i^{p^kj}}{j!}, 0\leq j <p\} \subseteq P_{s,i,k}$ is the
{\em iterative $\ad \, \der_i^{[p^k]}$-descent of maximal length}
$p$. Applying repeatedly Corollary \ref{1pw2Dec06}, we have the
projection
$$ \phi_{n+i}:= \phi_{n+i, k_i-1}\cdots \phi_{n+i, 1}\phi_{n+i,
0}: P_s\ra P_s, \;\; \phi_{n+i, k}:= \sum_{j=0}^{p-1} (-1)^j
\frac{x_i^{p^kj}}{j!} (\ad \, \der_i^{[p^k]})^j,$$ onto the
subalgebra $P_{s,i,k_i}$ of $P_s=\oplus_{j=0}^{p^{k_i}-1}
P_{s,i,k_i}x_i^j$.  The maps $\phi_{n+1}, \ldots , \phi_{2n}$ {\em
commute} and their product
$$ \prod_{i=1}^n \phi_{n+i}: P_s\ra P_s$$
is the projection onto $Z$ where $P_s=\oplus_{\alpha \in \CN}
Zx^\alpha$.

In order to avoid a clash of notations, let
$$ X_1:= x_1^{p^{k_1}}, \ldots , X_n:= x_n^{p^{k_n}},
X_{n+1}=x_{n+1}, \ldots , X_s = x_s.$$ Then  $Z= K[X_1, \ldots ,
X_s]$ is a polynomial algebra in $s$ variables. For each variable
$X_\nu$, let $\{ \der_\nu^{[j]}, j\geq 0\}\subseteq \CD
(Z)\subseteq {\rm End}_K(Z)$ be the corresponding higher
derivations (with respect to the variable $X_\nu $).

Let $\phi_Z$ be the map from (\ref{pphis}) in the case of the
polynomial algebra $Z=K[X_1, \ldots , X_s]$. Then the maps
%\marginpar{Tphi1}
\begin{equation}\label{Tphi1}
\phi_Z\phi_s\cdots\phi_{n+1} \phi_n\cdots \phi_1 : T\ra T, \;\;
a=\sum_{\alpha , \beta \in \CN , \g \in \Ns} \l_{ \beta \alpha\g }
\derb x^\alpha X^\g \mapsto \l_0, \;\; (\l_{\beta \alpha \g }\in
K),
\end{equation}
%\marginpar{Tpsi1}
\begin{equation}\label{Tpsi1}
\phi_Z \phi_s\cdots\phi_{n+1} \psi_n\cdots \psi_1 : T\ra T, \;\;
a=\sum_{\alpha , \beta \in \CN , \g \in \Ns} \l_{\alpha \beta \g }
x^\alpha\derb  X^\g \mapsto \l_0, \;\; (\l_{\alpha \beta \g }\in
K),
\end{equation}
are {\em projections} onto the field $K$,
$$ T= \bigoplus_{\alpha , \beta \in \CN , \g \in \Ns} Kx^\alpha\derb  X^\g
= \bigoplus_{\alpha , \beta \in \CN , \g \in \Ns} K\derb  x^\alpha
X^\g.$$ The maps %\marginpar{1Tphi1}
\begin{equation}\label{1Tphi1}
\phi :=\phi_s\cdots\phi_{n+1} \phi_n\cdots \phi_1 : T\ra T, \;\;
a=\sum_{\alpha , \beta \in \CN } a_{ \beta \alpha } \derb x^\alpha
\mapsto a_0, \;\; (a_{\beta \alpha  }\in Z),
\end{equation}
%\marginpar{1Tpsi1}
\begin{equation}\label{1Tpsi1}
\psi :=  \phi_s\cdots\phi_{n+1} \psi_n\cdots \psi_1 : T\ra T, \;\;
a=\sum_{\alpha , \beta \in \CN } a_{\alpha \beta } x^\alpha\derb
\mapsto a_0, \;\; (a_{\alpha \beta}\in Z),
\end{equation}
are {\em projections} onto the centre Z (see (\ref{TZab})).

\begin{theorem}\label{12Apr06}%\marginpar{12Apr06}
For any $a \in T$,
\begin{eqnarray*}
 a&=& \sum_{\alpha , \beta \in \CN , \g \in \Ns}(-1)^{|\alpha |}\phi_Z (\der^{[\g
 ]}(\phi ( \d^{\beta }  (a)\dera ))) \derb  x^\alpha  X^\g \\
 &=&\sum_{\alpha , \beta \in \CN , \g \in \Ns}\phi_Z (\der^{[\g ]}
 (\psi (
\dera \d^\beta  (a)))) x^\alpha \derb    X^\g ,
\end{eqnarray*}
where
 $\d^\beta:=
  \prod_{j=1}^n (-\ad \,
  x_j)^{\beta_j}$.
\end{theorem}

{\it Proof}. If $a= \sum_{\alpha , \beta \in \CN }  \l
_{\beta\alpha }\derb  x^\alpha  =\sum_{\alpha , \beta \in \CN } \l
_{ \alpha \beta }'x^\alpha \derb \in T$ where $\l _{
\beta\alpha},\l _{ \alpha \beta }'\in Z$, then it suffices to
prove that $(-1)^{|\alpha |}\phi ( \d^{\beta }  (a)\dera ))=
\l_{\beta \alpha }$ and $ \psi ( \dera \d^\beta  (a))= \l_{\alpha
\beta}'$ since, for any $z=\sum_{\g \in \Ns } z_\g X^\g $, $ z_\g
\in K$, we have $ \sum_{\g \in \Ns} \phi_Z (\der^{[\g ]}(z))=
z_0$. It suffices to prove the statements for $n=1$ (for an
arbitrary $n$, repeat the arguments below $n$ times). So, let
$n=1$. Recall that $\der_1^{[0]}:=1$ and $\der_1^{[s]}:=0$ for all
{\em negative} integers $s$.
\begin{eqnarray*}
 (-1)^i\phi (\d^j_1(a)\der_1^{[i]})&=& (-1)^i\phi ((-\ad \, x_1)^j
 (\sum_{\alpha, \beta \in \CN}\l_{\beta \alpha }\der_1^{[\beta ]}x_1^\alpha )\der_1^{[i]})
 = (-1)^i\phi (\sum_{\alpha, \beta \in \CN, \beta \geq j}\l_{\beta \alpha
 }\der_1^{[\beta-j ]}x_1^\alpha \der_1^{[i]}) \\
 &=&(-1)^i\phi (\sum_{\alpha \in \CN}\l_{j \alpha }x_1^\alpha
 \der_1^{[i]})=(-1)^i\phi (\sum_{\alpha \in \CN}\l_{j \alpha }\sum_{k=0}^\alpha {\alpha \choose k} (-1)^k
 \der_1^{[i-k]}x_1^{\alpha-k}) \\
 &=&(-1)^i(-1)^i\l_{ji}=\l_{ji},
\end{eqnarray*}
\begin{eqnarray*}
 \psi (\der_1^{[i]}\d^j_1(a))&=& \psi (\der_1^{[i]}(-\ad\,
 x_1)^j(a))= \psi  (\der_1^{[i]}(-\ad \, x_1)^j
 (\sum_{\alpha, \beta \in \CN }\l_{ \alpha \beta}'x_1^\alpha\der_1^{[\beta ]} ))\\
&=& \psi  (\der_1^{[i]}
 \sum_{\alpha, \beta \in \CN , \beta \geq j}\l_{ \alpha \beta}'x_1^\alpha\der_1^{[\beta -j]}
 )= \psi  (
 \sum_{\alpha \in \CN }\l_{ \alpha j}'\der_1^{[i]}x_1^\alpha
 )\\
 &=&\psi ( \sum_{\alpha \in \CN }\l_{ \alpha
j}' \sum_{k=0}^\alpha {\alpha \choose k}x_1^{\alpha -k}
\der_1^{[i-k]}) =\l_{ij}'. \;\; \Box
\end{eqnarray*}

Let $\s\in \Aut_K(T)$ and let $x_i':=\s (x_i)$,
${\der_i'}^{[p^{s_i}]}:= \s(\der_i^{[p^{s_i}]})$, $i=1,\ldots ,
n$, $0\leq s_i<k_i$,  and $X_j':= \s (X_j)$, $j=1, \ldots , s$.
Then $Z=K[X_1', \ldots X_s']$.  Let $\phi_{Z, \s }$, $\phi_\s $
and $\psi_\s$ be the maps (\ref{pphis}),  (\ref{1Tphi1}) and
(\ref{1Tpsi1}) for the new choice of the canonical generators:
for, we have to put $(')$ everywhere.

\begin{theorem}\label{iT12Apr06}%\marginpar{iT12Apr06}
({\rm The
inversion formula for} $\s \in \Aut_K(T)$)  For any $\s \in
\Aut_K(T)$ and  $a \in T$,
\begin{eqnarray*}
\s^{-1} ( a) &=& \sum_{\alpha , \beta \in \CN , \g \in
\Ns}(-1)^{|\alpha |}\phi_{Z, \s } ({\der'}^{[\g
 ]}(\phi_\s  ( {\d'}^{\beta }  (a){\der'}^{[\alpha]} ))) \derb  x^\alpha  X^\g \\
 &=&\sum_{\alpha , \beta \in \CN , \g \in \Ns}\phi_{Z, \s } ({\der'}^{[\g ]}
 (\psi_\s  (
{\der'}^{[\alpha ]} {\d'}^\beta  (a)))) x^\alpha \derb    X^\g ,
\end{eqnarray*}
where  ${\der'}^{[\g ]}:=
  \prod_{i=1}^s{\der'}_i^{[\g_i]}$ and
 ${\d'}^\beta:=
  \prod_{j=1}^n (-\ad \,
  {x'}_j)^{\beta_j}$.

\end{theorem}

{\it Proof}. By Theorem \ref{12Apr06},
\begin{eqnarray*}
\s^{-1} ( a) &=& \sum_{\alpha , \beta \in \CN , \g \in
\Ns}(-1)^{|\alpha |}\phi_{Z, \s } ({\der'}^{[\g
 ]}(\phi_\s  ( {\d'}^{\beta }  (a){\der'}^{[\alpha]} ))) {\der'}^{[\beta ]} {x'}^\alpha  {X'}^\g \\
 &=&\sum_{\alpha , \beta \in \CN , \g \in \Ns}\phi_{Z, \s } ({\der'}^{[\g ]}
 (\psi_\s  (
{\der'}^{[\alpha ]} {\d'}^\beta   (a)))) {x'}^\alpha
{\der'}^{[\beta ]} {X'}^\g ,
\end{eqnarray*}
 then applying $\s^{-1}$ proves the result.  $\Box$

Department of Pure Mathematics

University of Sheffield

Hicks Building

Sheffield S3 7RH

UK

email: v.bavula@sheffield.ac.uk

\end{document}